# EXISTENCE OF THE ZERO RANGE PROCESS AND A DEPOSITION MODEL WITH SUPERLINEAR GROWTH RATES

By M. Balázs,[1] F. Rassoul-Agha,[2] T. Seppäläinen[3]
and S. Sethuraman[4]

*University of Wisconsin–Madison, University of Utah, University of Wisconsin–Madison and Iowa State University*

We give a construction of the zero range and bricklayers' processes in the totally asymmetric, attractive case. The novelty is that we allow jump rates to grow exponentially. Earlier constructions have permitted at most linearly growing rates. We also show the invariance and extremality of a natural family of i.i.d. product measures indexed by particle density. Extremality is proved with an approach that is simpler than existing ergodicity proofs.

**1. Introduction.** Balázs [2, 5] introduced a stochastic deposition model named the *bricklayers' process*. This process represents a wall formed by adjacent columns of bricks, where each column grows with a rate that depends on the relative heights of the neighboring columns. There are other representations that show close connection to interacting particle systems; the details will be given later. The process with exponentially growing rates is of special interest because this case possesses shock-like product-form invariant distributions, as discovered in [5]. In the present paper we address the rigorous construction and the ergodicity of the process which were assumed but left open in the earlier work.

In the area of interacting particle systems there are two main situations where construction methods are available. In the first situation the total

Received November 2005; revised October 2006.
[1]Supported in part by the Hungarian Scientific Research Fund (OTKA) Grants TS49835 and T037685 and NSF Grant DMS-05-03650.
[2]Supported in part by NSF Grant DMS-05-05030.
[3]Supported in part by NSF Grant DMS-04-02231.
[4]Supported in part by NSF Grant DMS-05-04193 and NSA Grant H982300510041.

rate of change at a site is bounded. As described in [8], in this case the construction can be based on functional analytic properties of the infinitesimal generator and the Hille–Yosida theorem. This way existence of dynamics can be proved for stochastic Ising models, the voter model, contact processes, and simple- and $K$-exclusion processes.

In the other situation the rate is bounded above by a linear function of the local state space. The most-studied example is the zero range process. The state of this process is a configuration of occupation numbers $(\omega_i)$. When the number of particles at site $i$ is $\omega_i$, one particle is moved to another site at rate $r(\omega_i)$. Assume the Lipschitz condition $|r(k+1)-r(k)| \leq K$ for all $k \geq 0$, with some constant $0 < K < \infty$. This of course implies that $r(k)$ grows at most linearly. Under this condition, the zero range process can be compared to a multitype branching process to get stochastic bounds on the evolution. This way Andjel [1] constructed the zero range process, generalizing earlier work of Liggett [7].

The approach of Andjel can be extended to more complicated systems, but (sub)linearity is still an essential condition of the proof of existence. Booth [6] and Quant [10] apply these ideas to the bricklayers' process. But neither approach mentioned above works for the bricklayers' process with exponential rates which is the process studied in [2] and [5].

Let us also mention convexity, which naturally appears together with superlinearity in certain rate functions, such as polynomials or exponentials. Convexity plays an important role in hydrodynamical limits and arguments that involve second class particles. For examples, see [3] and [11].

In the present paper we consider the bricklayers' process with unbounded jump rates and relax the linear growth condition to an *exponential* one. We use attractivity and couplings to construct the infinite system as a limit of systems on finitely many sites. Then we develop some properties of this process. We associate with it the natural generator suggested by an informal description of the dynamics, and we prove ergodicity and extremality of a family of i.i.d. probability measures on the configurations.

For the case of linearly bounded rates, [1] introduced a norm and a state space precisely tailored to the rates. The evolution of the norm could be controlled suitably, and then a semigroup constructed on the space of bounded Lipschitz functions on the state space. By contrast, in the present case we did not find a suitable norm, so we take as our state space $\widetilde{\Omega}$ the space of configurations whose asymptotic slopes are Cesàro-bounded in a sense to be made precise. The attractive feature is that a single state space works for all choices of rates. The drawback is that we have no useful topology on $\widetilde{\Omega}$, and we regard it only as a Borel subset of the product space $\Omega = \mathbb{Z}^{\mathbb{Z}}$ of all configurations. As a consequence, the analytic properties we prove for the semigroup and the generator fall short of what can be achieved for bounded and linearly bounded rates.

Two basic tools we utilize throughout the paper are the following: (1) a *conditional coupling*, where we condition two coupled processes on a favorable event at time zero, and (2) bounds on the probability that a large block of adjacent sites all experience a jump in a given time interval. The conditional coupling is used to bound processes started from arbitrary initial states with processes started from conveniently chosen initial distributions. The bounds on propagation of jumps allow us to bound the motion of second class particles. Such bounds are needed for analytic properties of the semigroup and for the proof of stationarity and ergodicity.

Our paper treats explicitly only the bricklayers' process. The arguments work also for the totally asymmetric zero range process. Hence, we obtain a construction of this process for monotone rates that are bounded by an exponential. We emphasize that only the totally asymmetric version of the zero range process is covered by this paper.

**2. The process and the results.** We now define bricklayers' process precisely, address briefly the application of this paper to the somewhat simpler zero-range case, and then describe the results. Subsequent sections of the paper develop the proofs, beginning with definitions of processes on finitely many sites.

2.1. *Formal description of the bricklayers' process.* Before we define a restricted state space for the process, the states are simply elements of the product space $\Omega = \mathbb{Z}^{\mathbb{Z}}$. In other words, a state is defined by attaching to each site $i \in \mathbb{Z}$ an integer $\omega_i \in \mathbb{Z}$. Let $r : \mathbb{Z} \to \mathbb{R}_+$ be a strictly positive function. Given a configuration

$$\underline{\omega} = \{\omega_i \in \mathbb{Z} : i \in \mathbb{Z}\} \in \Omega,$$

define another configuration $\underline{\omega}^{(i,i+1)}$ by

$$(\underline{\omega}^{(i,i+1)})_j = \begin{cases} \omega_j, & \text{for } j \neq i, i+1, \\ \omega_i - 1, & \text{for } j = i, \\ \omega_{i+1} + 1, & \text{for } j = i+1. \end{cases}$$

The informal description of the evolution is this: conditioned on the present state $\underline{\omega}$, the jump $\underline{\omega} \to \underline{\omega}^{(i,i+1)}$ happens independently for each site $i$ with rate $r(\omega_i) + r(-\omega_{i+1})$. A pictorial interpretation of this jump follows below.

We assume four conditions on the rate function $r$:

- $r$ is a strictly increasing function. A nondecreasing $r$ makes the process *attractive*, which is useful for coupling arguments. A strictly increasing $r$ is required only for the ergodicity arguments in Section 7. The remainder of the paper needs only a nondecreasing $r$.

- To have product invariant distributions, we require

(2.1) $$r(z) \cdot r(1-z) = 1$$

  for each $z \in \mathbb{Z}$. Note that now the rates for positive $z$'s determine the rates for negative $z$'s.

REMARK 2.1. In the totally asymmetric zero range process $r(z) = 0$ for all $z \leq 0$. Hence, the right-hand side of (2.1) is zero. In Section 2.2 below we introduce briefly the zero range process.

- The condition

(2.2) $$\lim_{z \to \infty} r(z) = \infty$$

  is needed for fast enough tail decay of the invariant distributions.
- Finally, the *exponential* growth condition that is the key point of the paper: there is a constant $\beta > 0$ such that

(2.3) $$r(z) < e^{\beta z} \qquad \text{for all } z > 0.$$

The informal description of the dynamics suggests the following formal infinitesimal generator that acts on test functions $\varphi : \Omega \to \mathbb{R}$:

$$L\varphi(\underline{\omega}) = \sum_{i \in \mathbb{Z}} [r(\omega_i) + r(-\omega_{i+1})] \cdot [\varphi(\underline{\omega}^{(i,i+1)}) - \varphi(\underline{\omega})].$$

$L\varphi$ is well defined for bounded cylinder functions.

Next we describe a picture associated with the dynamics that motivates the name "bricklayers' process." Imagine a wall consisting of side-by-side columns of bricks. Each column occupies the width between two consecutive integer sites. The variable $\omega_i$ is the negative discrete gradient of the wall at site $i$, in other words, the difference between the heights of the columns to the left of $i$ and to the right of $i$. The jump $\underline{\omega} \to \underline{\omega}^{(i,i+1)}$ corresponds to the growth of the column over $[i, i+1]$; see Figure 1. Thus, the growth rate for each column consists of two additive parts depending on $\omega$ of the left-hand side and on $\omega$ of the right-hand side, respectively. The process can be represented by an infinite row of conditionally independent bricklayers. Each bricklayer occupies his own site $i$ and lays bricks to his right with rate $r(\omega_i)$ and to his left with rate $r(-\omega_i)$.

It is convenient to work with the heights because they increase monotonically, so let us introduce notation for them. The height of the column between sites $i$ and $i+1$ is denoted by $h_i$. The connection with the increment variables is

(2.4) $$h_i = \begin{cases} h_0 - \sum_{j=1}^{i} \omega_j, & \text{for } i > 0, \\ h_0 + \sum_{j=i+1}^{0} \omega_j, & \text{for } i < 0. \end{cases}$$

The jump $\underline{\omega} \to \underline{\omega}^{(i,i+1)}$ is associated to the jump $h_i \to h_i + 1$. Once the processes are well defined, the variables $\omega_i(t)$ and $h_i(t)$ satisfy (2.4) for all later times $t$ also. When we are only interested in the increment process $\underline{\omega}(t)$, it is convenient to normalize the heights initially so that $h_0(0) = 0$. But, in general, arbitrary height configurations $\underline{h} = (h_i)$ are permitted. The attractivity, which follows from the monotonicity of the function $r$, has an especially natural meaning for the height process: the higher the neighbors of a column, the faster this column grows.

We finish this subsection with a definition of the measures that will be shown invariant. Define
$$r(0)! := 1, \qquad r(n)! := \prod_{y=1}^{n} r(y), \qquad n \in \mathbb{N}.$$

For a real parameter $\theta$, define
$$Z(\theta) := \sum_{z=-\infty}^{\infty} \frac{e^{\theta z}}{r(|z|)!}.$$

The series converges by (2.2). Hence, we can define the probability measure

(2.5) $$\mu^{(\theta)}(z) := \frac{1}{Z(\theta)} \cdot \frac{e^{\theta z}}{r(|z|)!}$$

on the space $\mathbb{Z}$ of integers. On the product space $\Omega = \mathbb{Z}^{\mathbb{Z}}$ define the i.i.d. product measure $\underline{\mu}^{(\theta)}$ with marginals

(2.6) $$\underline{\mu}^{(\theta)}\{\underline{\omega} : \omega_i = z\} = \mu^{(\theta)}(z).$$

Expectation w.r.t. $\underline{\mu}^{(\theta)}$ or $\mu^{(\theta)}$ will be denoted by $\mathbf{E}^{(\theta)}$. A formal computation with the generator $L$ suggests that the product measure $\underline{\mu}^{(\theta)}$ is stationary for the process. We prove this in Section 7. Some properties of these measures have been relegated to the Appendix.

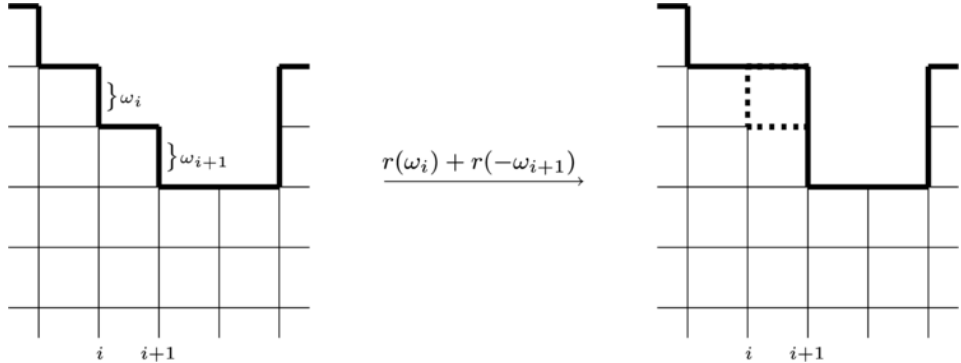

FIG. 1. *A possible move.*

2.2. *The totally asymmetric zero range process.* The increment process $\underline{\omega}(t)$ of the bricklayers' process can be given a particle description. Suppose $\omega_i$ represents the number of particles at site $i$. Then we have a nearest-neighbor process, where particles ($\omega_i > 0$) jump to the right, antiparticles ($\omega_i < 0$) jump to the left, and the two types annihilate each other upon meeting.

The totally asymmetric zero range process is quite similar to the bricklayers' process, but without the antiparticles. So now $\Omega = (\mathbb{Z}^+)^{\mathbb{Z}}$, that is, $\omega_i \geq 0$ for all $i \in \mathbb{Z}$. The formal infinitesimal generator is

$$(L\varphi)(\underline{\omega}) = \sum_{i \in \mathbb{Z}} r(\omega_i) \cdot [\varphi(\underline{\omega}^{(i,i+1)}) - \varphi(\underline{\omega})].$$

Instead of (2.1), we require $r(0) = 0$. Then $\omega_i \geq 0$ can never be violated. The measure $\mu^{(\theta)}$, now a product of the marginal measures (2.5) restricted to $\mathbb{Z}^+$, is (formally) stationary for the process. Throughout the paper we show the arguments only for the bricklayers' process. Often these arguments will work for the zero range process with little or no change. One simply neglects all terms of the type $r(-\omega_i)$, $r(-\zeta_i)$, $r(-\xi_i)$, $r(-\eta_i)$, $e^{-\theta}$, $e^{-\theta_1}$, $e^{-\theta_2}$. We shall remark on points that need different treatment for the zero range process.

2.3. *Results for bricklayers' process.* We now discuss the rigorous results. First we introduce a process over a finite number of sites. Fix integers $\ell < \mathfrak{r}$, and define the infinitesimal generator $L^{[\ell, \mathfrak{r}]}$ acting on functions of $\underline{\omega}$:

$$(2.7) \qquad L^{[\ell, \mathfrak{r}]} \varphi(\underline{\omega}) = \sum_{i=\ell}^{\mathfrak{r}-1} [r(\omega_i) + r(-\omega_{i+1})] \cdot [\varphi(\underline{\omega}^{(i,i+1)}) - \varphi(\underline{\omega})].$$

This generator defines a countable state space process evolving over the sites $\ell, \ldots, \mathfrak{r}$: the jump $\underline{\omega} \to \underline{\omega}^{(i,i+1)}$ happens with rate $r(\omega_i) + r(-\omega_{i+1})$, independently for different sites $i$, but only for $\ell \leq i \leq \mathfrak{r} - 1$. Columns outside the interval $[\ell, \mathfrak{r}]$ are frozen for all time. We show in Section 3.1 that this process, denoted by $\underline{\omega}^{[\ell, \mathfrak{r}]}(t)$, is well defined for any initial state $\underline{\omega} \in \Omega$. We call this process the $[\ell, \mathfrak{r}]$-*monotone process*, and we justify this name in Section 3.4. The bricklayers' process will be produced as the infinite volume limit of the $[\ell, \mathfrak{r}]$-monotone process as $\ell \to -\infty$ and $\mathfrak{r} \to \infty$:

THEOREM 2.2. *Fix $\underline{\omega} \in \Omega$. Consider the family $\{\underline{\omega}^{[\ell, \mathfrak{r}]}(t) : -\infty < \ell < \mathfrak{r} < \infty\}$ of all $[\ell, \mathfrak{r}]$-monotone processes with the common initial state $\underline{\omega} \in \Omega$. The height of column $i$ of $\underline{\omega}^{[\ell, \mathfrak{r}]}(t)$ at time $t$ is denoted by $h_i^{[\ell, \mathfrak{r}]}(t)$, with the initial normalization $h_0^{[\ell, \mathfrak{r}]}(0) = 0$. Then there is a coupling such that, for all $\mathcal{L} \leq \ell < 0 < \mathfrak{r} \leq \mathcal{R}$, $h_i^{[\ell, \mathfrak{r}]}(t) \leq h_i^{[\mathcal{L}, \mathcal{R}]}(t)$ for all $i \in \mathbb{Z}$ and $0 \leq t < \infty$. Consequently, the monotone nondecreasing limits*

$$h_i(t) := \lim_{[\ell, \mathfrak{r}] \nearrow \mathbb{Z}} h_i^{[\ell, \mathfrak{r}]}(t) = \sup_{[\ell, \mathfrak{r}] \subseteq \mathbb{Z}} h_i^{[\ell, \mathfrak{r}]}(t)$$

*exist a.s. for all $i \in \mathbb{Z}$ and $0 \leq t < \infty$. Since we do not restrict the initial state $\underline{\omega}$, $h_i(t) = \infty$ is possible at this point.*

This theorem is basic for everything that follows: *throughout our paper, unless otherwise stated, the height process $\underline{h}(\cdot) = (h_i(t) : i \in \mathbb{Z}, t \geq 0)$ is constructed by the limit in the above theorem.* The increment process $\underline{\omega}(t)$ is defined in terms of the height process by $\omega_i(t) := h_{i-1}(t) - h_i(t)$ whenever the heights are finite.

Naturally, next we seek to restrict the initial states $\underline{\omega}$ in order to get well-defined, finite-valued, infinite-volume processes. First we consider initial states distributed according to certain special measures. Recall the definition (2.6) of $\underline{\mu}^{(\theta)}$. Let us say that a probability measure $\underline{\pi}$ on $\Omega$ is a *good measure with parameters* $\theta_1 < \theta_2$ if the stochastic dominations $\underline{\mu}^{(\theta_2)} \geq \underline{\pi} \geq \underline{\mu}^{(\theta_1)}$ hold. Equivalently, there is a coupling of three random configurations

$$\underline{\eta} \sim \underline{\mu}^{(\theta_1)}, \quad \underline{\zeta} \sim \underline{\pi} \quad \text{and} \quad \underline{\xi} \sim \underline{\mu}^{(\theta_2)}$$

such that

$$\eta_i \leq \zeta_i \leq \xi_i \quad \text{almost surely for all } i \in \mathbb{Z}.$$

If $\underline{\pi}$ is a product of marginals $\pi_i$ on $\mathbb{Z}$, then this is equivalent to the corresponding stochastic domination for the marginals at each site $i$. An obvious example of a good measure with parameters $\theta_1$ and $\theta_2$ is $\underline{\mu}^{(\theta)}$ for $\theta \in [\theta_1, \theta_2]$ (Lemma A.2 in the Appendix).

THEOREM 2.3. *Let the initial state $\underline{\omega}$ be distributed according to a good measure $\underline{\pi}$ with parameters $\theta_1 < \theta_2$, and consider the processes $\underline{\omega}^{[\ell, \mathrm{r}]}(t)$ started from the common initial configuration $\underline{\omega}$. Then the height process constructed in Theorem 2.2 satisfies*

$$\mathbf{E}[h_i(t) - h_i(0)] \leq t \cdot (e^{\theta_2} + e^{-\theta_1}) \quad \text{for all } i \in \mathbb{Z} \text{ and } t \geq 0.$$

Once we have control over processes with good random initial states, we can turn to control processes with suitably restricted deterministic initial states. Define

$$(2.8) \quad K_{\underline{\omega}} := \max\left\{ \limsup_{i \to -\infty} \frac{1}{|i|} \sum_{j=i+1}^{0} |\omega_j|, \limsup_{i \to \infty} \frac{1}{i} \sum_{j=1}^{i} |\omega_j| \right\} \in [0, \infty].$$

Our state space for the increment process will be

$$(2.9) \quad \widetilde{\Omega} := \{\underline{\omega} \in \Omega : K_{\underline{\omega}} < \infty\}.$$

We could call this the space of configurations with Cesàro-bounded slopes. We endow $\widetilde{\Omega}$ with the product topology and product $\sigma$-algebra it inherits from the product space $\Omega$. It is simple but important that $\underline{\mu}^{(\theta)}(\widetilde{\Omega}) = 1$ for all $\theta \in \mathbb{R}$. Let us denote expectation (probability) for processes $\underline{\omega}^{[\ell, \mathrm{r}]}(t)$ or $\underline{\omega}(t)$ started from initial state $\underline{\omega}$ by $\mathbf{E}^{\underline{\omega}}$ ($\mathbf{P}^{\underline{\omega}}$, resp.).

THEOREM 2.4. *Suppose the initial state $\underline{\omega} \in \widetilde{\Omega}$. Then there are finite constants $A^{\underline{\omega}}$ and $B^{\underline{\omega}}$ that depend on the initial state $\underline{\omega}$ such that the height process constructed in Theorem 2.2 satisfies*

$$\mathbf{E}^{\underline{\omega}}[h_i(t)] \leq A^{\underline{\omega}} t + B^{\underline{\omega}} |i| \qquad \text{for all } 0 \leq t < \infty \text{ and all } i \in \mathbb{Z}.$$

In the settings of the above two theorems we have well-defined, finite-valued, infinite-volume processes $\underline{h}(t) = (h_i(t))_{i \in \mathbb{Z}}$ and $\underline{\omega}(t) = (\omega_i(t))_{i \in \mathbb{Z}}$ for all time $t \in [0, \infty)$. The next theorem justifies calling $\widetilde{\Omega}$ of (2.9) the state space. Let $D_{\widetilde{\Omega}}$ denote the space of cadlag functions from $[0, \infty)$ into $\widetilde{\Omega}$, endowed with the $\sigma$-algebra generated by the coordinates.

THEOREM 2.5. *Suppose the initial increments satisfy $\underline{\omega} \in \widetilde{\Omega}$. Then $\mathbf{P}^{\underline{\omega}}$-almost surely $\underline{\omega}(\cdot)$ is a Markov process with paths in $D_{\widetilde{\Omega}}$. Furthermore, we have this uniformity for the asymptotic slope: there is a finite real $\widehat{K}_{\underline{\omega}}$ such that*

$$\mathbf{P}^{\underline{\omega}} \left\{ \sup_{n \in \mathbb{Z}_+} \limsup_{i \to \infty} \sup_{t \in [n, n+1]} \frac{1}{i} \sum_{j=-i}^{i} |\omega_j(t)| \leq \widehat{K}_{\underline{\omega}} \right\} = 1.$$

In other words, we cannot show that the constant $K_{\underline{\omega}}$ in (2.8) of the initial state is preserved by the evolution, but another larger constant $\widehat{K}_{\underline{\omega}}$ is preserved. Elements of $D_{\widetilde{\Omega}}$ will be denoted by $\underline{\omega}(\cdot)$. There is no need to distinguish notationally between the process and elements of its path space. By similar abuse of notation, $\mathbf{P}^{\underline{\omega}}$ denotes also the path measure on $D_{\widetilde{\Omega}}$ for the process $\underline{\omega}(\cdot)$ with initial state $\underline{\omega}$. Along the way, it will be verified that $\mathbf{P}^{\underline{\omega}}$ is measurable as a function of $\underline{\omega}$, in the sense that an expectation

$$\mathbf{E}^{\underline{\omega}}[F] = \int_{D_{\widetilde{\Omega}}} F(\underline{\omega}(\cdot)) \mathbf{P}^{\underline{\omega}}(d\underline{\omega}(\cdot))$$

of a bounded measurable function $F$ on $D_{\widetilde{\Omega}}$ is a measurable function of the initial state $\underline{\omega}$.

A semigroup of linear contractions acting on the space of bounded measurable functions on $\widetilde{\Omega}$ is defined by

$$S(t) f(\underline{\omega}) = \mathbf{E}^{\underline{\omega}}[f(\underline{\omega}(t))].$$

As usual, the semigroup property means that $S(0) = I$ and $S(s+t) = S(s)S(t)$. These operators are contractions if functions are normed by the sup-norm. The approximation in Theorem 2.2 will have the property that $\omega_i(t) = \omega_i^{[\ell, \mathfrak{r}]}(t)$ for large enough $-\ell, \mathfrak{r}$ when $(i, t)$ are restricted to a compact set and the initial state is from $\widetilde{\Omega}$. This is used throughout the proofs.

The height process $\underline{h}(\cdot)$ is also a Markov process with paths in a $D$-space. In fact, in the construction $\underline{h}(\cdot)$ comes first, so we shall have more to say about the height process in later sections.

Theorems 2.2–2.4 above are simple consequences of a series of coupling arguments contained in Sections 3 and 4. Theorem 2.5 requires slightly more involved estimates and is proved in Section 5.

Section 6 turns to analytic properties of the semigroup constructed above. Here is a summary of the results. Throughout, $\varphi$ is a bounded cylinder function on $\widetilde{\Omega}$, in other words, a bounded function of some finite set $(\omega_i)_{-a \leq i \leq a}$ of variables. Also, the starting state $\underline{\omega}$ lies in $\widetilde{\Omega}$.

THEOREM 2.6. *The function $F(t,\underline{\omega}(\cdot)) = |L\varphi(\underline{\omega}(t))|$ defined on $[0,T] \times D_{\widetilde{\Omega}}$ is integrable under the measure $dt \otimes \mathbf{P}^{\underline{\omega}}$, for any $0 < T < \infty$. Thus, $S(t)L\varphi(\underline{\omega}) = \mathbf{E}^{\underline{\omega}}[L\varphi(\underline{\omega}(t))]$ is well defined for almost every $t > 0$ and integrable over any bounded time interval $[0,T]$. For all $t > 0$, we have the integrated forward equation*

$$S(t)\varphi(\underline{\omega}) = \varphi(\underline{\omega}) + \int_0^t S(s)L\varphi(\underline{\omega})\,ds.$$

The integrated backward equation and the differential equation are only proved for a short time that depends on the initial state.

THEOREM 2.7. *There exists $T^{\underline{\omega}} > 0$ such that, for any $t \in (0, T^{\underline{\omega}}]$, the following statements hold:*

(i) *$LS(t)\varphi(\underline{\omega})$ is well defined and the integrated backward equation holds:*

$$S(t)\varphi(\underline{\omega}) = \varphi(\underline{\omega}) + \int_0^t LS(s)\varphi(\underline{\omega})\,ds.$$

(ii) *$S(t)L\varphi(\underline{\omega})$ is well defined, and there is a continuous derivative*

$$\frac{d}{dt}S(t)\varphi(\underline{\omega}) = S(t)L\varphi(\underline{\omega}) = LS(t)\varphi(\underline{\omega}).$$

*In particular, $(d/dt)S(t)\varphi(\underline{\omega})|_{t=0} = L\varphi(\underline{\omega})$.*

The approximation $\underline{\omega}^{[\ell,\mathrm{r}]}(t) \to \underline{\omega}(t)$ and Theorems 2.6 and 2.7 are valid over the entire state space rather than only locally, if the dependence on the spatial tails is weak enough. For example, if we define a norm on $\widetilde{\Omega}$ by $\|\underline{\omega}\| = \sum_{i \in \mathbb{Z}} e^{-i^2}|\omega_i|$, then $\underline{\omega}^{[\ell,\mathrm{r}]}(t) \to \underline{\omega}(t)$ is valid in this norm, and so are Theorems 2.6 and 2.7 for bounded functions $\varphi : \widetilde{\Omega} \to \mathbb{R}$ that are Lipschitz continuous in this norm. We omit the details of this generalization since we have no further use for it. (These details can be found in the first version of this paper on ArXiv.)

We conclude the results with the invariance and ergodicity. Andjel [1] and Sethuraman [14] characterized the invariant and extremal invariant measures of the zero range process. Our present goal is only to prove the invariance and extremality of the product measures $\underline{\mu}^{(\theta)}$ defined in (2.6) for the process $\underline{\omega}(\cdot)$ constructed by the limit in Theorem 2.2. Let

$$\mathbf{P}^{\underline{\mu}^{(\theta)}} = \int_{\widetilde{\Omega}} \mathbf{P}^{\underline{\omega}} \underline{\mu}^{(\theta)}(d\underline{\omega})$$

denote the path measure on the space $D_{\widetilde{\Omega}}$ with initial distribution $\underline{\mu}^{(\theta)}$. As usual, stationarity means the invariance $\mathbf{P}^{\underline{\mu}^{(\theta)}} \circ \sigma_t^{-1} = \mathbf{P}^{\underline{\mu}^{(\theta)}}$ under time shifts defined by $(\sigma_t \underline{\omega})(s) = \underline{\omega}(s+t)$. Ergodicity means that the shift-invariant sets are trivial: $\mathbf{P}^{\underline{\mu}^{(\theta)}}(A) = 0$ or 1 for events $A$ on $D_{\widetilde{\Omega}}$ that satisfy $\sigma_t^{-1} A = A$ for all $0 \le t < \infty$.

THEOREM 2.8. *For each $\theta \in \mathbb{R}$, the path measure $\mathbf{P}^{\underline{\mu}^{(\theta)}}$ is stationary and ergodic under time shifts.*

Of the various equivalent characterizations of ergodicity, we verify that harmonic functions in $L^2(\underline{\mu}^{(\theta)})$ are a.e. constant. This is done in Section 7 with a fairly simple approach. "Extremality" of $\underline{\mu}^{(\theta)}$ above refers to the fact that extremality of $\underline{\mu}^{(\theta)}$ among invariant distributions is one equivalent characterization of ergodicity. For more about this, we refer to [12] and [14].

2.4. *Open problems.* Our estimates are such that many bounds hold only for a time range that depends on the initial state. A better approach should be found. This should enable one to derive better analytic properties of the semigroup.

Since our columns only grow and never decrease, this model corresponds to a totally asymmetric particle system. (In particular, only the totally asymmetric case of a zero range process is treated.) It would be of interest to see if the arguments generalize to systems where bricks are both added and removed.

### 3. The process on a finite number of sites.

3.1. *The $[\ell, \mathfrak{r}]$-monotone process.* We begin by observing that the finite-volume process $\underline{\omega}^{[\ell, \mathfrak{r}]}(\cdot)$ generated by $L^{[\ell, \mathfrak{r}]}$ of (2.7) is well defined. Such a countable state space Markov process is well defined for all time if the process avoids explosion for each initial configuration $\underline{\omega}(0)$ (see Chapter 2 of [9]). We achieve this by showing that the corresponding height process $\underline{h}(t)$ is stochastically dominated. Since the heights never decrease, upper bounds suffice. Let

$$H(t) := \max_{\ell \leq j \leq \mathfrak{r}-1} h_j(t) \quad \text{and} \quad J := \{\ell \leq j \leq \mathfrak{r}-1 : h_j(t) = H(t)\}.$$

Then we have two possibilities:

(i) $H(t) \leq h_{\ell-1}(0) = h_{\ell-1}(t)$ or $H(t) \leq h_\mathfrak{r}(0) = h_\mathfrak{r}(t)$ (these are the heights of the closest columns to the origin which do not grow), or

(ii) for each $j \in J$, $\omega_j \leq 0$ and $\omega_{j+1} \geq 0$, since $h_j(t)$ is by definition maximal.

In the first case $H(t)$ is trivially dominated. In the second case, by monotonicity of the rates $r$, for each $j \in J$, the column between sites $j$ and $j+1$ has a growth rate dominated by $r(0) + r(0)$. This implies that whenever $H(t) \geq \max(h_{\ell-1}(0), h_\mathfrak{r}(0))$, it grows according to a continuous time jump process with rate dominated by $(\mathfrak{r} - \ell) \cdot (r(0) + r(0))$. This bound is good for each fixed $[\ell, \mathfrak{r}]$ but useless for the limit $[\ell, \mathfrak{r}] \nearrow \mathbb{Z}$.

3.2. *The $(\ell, \mathfrak{r}, \theta)$-process.* To develop equilibrium bounds, we introduce another finite-volume process whose rates are adjusted at the boundary. For integers $\ell < \mathfrak{r}$ and $\theta \in \mathbb{R}$, the generator is

$$\begin{aligned}
G^{(\ell,\mathfrak{r},\theta)}\varphi(\underline{\omega}) &= \sum_{i=\ell}^{\mathfrak{r}-1}[r(\omega_i) + r(-\omega_{i+1})] \cdot [\varphi(\underline{\omega}^{(i,i+1)}) - \varphi(\underline{\omega})] \\
&\quad + [e^\theta + r(-\omega_\ell)] \cdot [\varphi(\underline{\omega}^{(\ell-1,\ell)}) - \varphi(\underline{\omega})] \\
&\quad + [e^{-\theta} + r(\omega_\mathfrak{r})] \cdot [\varphi(\underline{\omega}^{(\mathfrak{r},\mathfrak{r}+1)}) - \varphi(\underline{\omega})].
\end{aligned} \quad (3.1)$$

In contrast with $L^{[\ell,\mathfrak{r}]}$, in $G^{(\ell,\mathfrak{r},\theta)}$ the bricklayer on site $\ell$ works also to his left, not only to his right, and similarly, the bricklayer at site $\mathfrak{r}$ also lays bricks to his right. Additionally, there are two "virtual" bricklayers. One at site $\ell - 1$ lays bricks on $[\ell - 1, \ell]$ at the mean equilibrium rate $e^\theta$ independently of the surrounding configuration. The other at site $\mathfrak{r} + 1$ lays bricks on $[\mathfrak{r}, \mathfrak{r}+1]$ at rate $e^{-\theta}$. These rates are the expected $r(\pm \omega_i)$ rates under $\mu^{(\theta)}$, hence, they can be called mean equilibrium rates.

By the argument of Section 3.1, the height of any column in the growing area is dominated by a Poisson process of rate $(\mathfrak{r} - \ell) \cdot (r(0) + r(0)) + 2r(0) + e^\theta + e^{-\theta}$ (we do not even need to consider the nongrowing columns, as this process is not sensitive to its surrounding environment). We call this process the $(\ell, \mathfrak{r}, \theta)$-process.

PROPOSITION 3.1. *The product-measure*

$$\underline{\mu}^{(\ell,\mathfrak{r},\theta)}(\underline{\omega}) := \prod_{i=\ell}^{\mathfrak{r}} \mu^{(\theta)}(\omega_i) \quad (3.2)$$

*is stationary for the $(\ell, \mathfrak{r}, \theta)$-process, where the marginals are of the form* (2.5).

PROOF. Since the growth rates only depend on $\omega_i$ for $\ell \leq i \leq \mathfrak{r}$, it is sufficient to show that, with the expectation w.r.t. $\underline{\mu}^{(\ell,\mathfrak{r},\theta)}$ of (3.2),

$$\mathbf{E}^{(\ell,\mathfrak{r},\theta)}((G^{(\ell,\mathfrak{r},\theta)}\varphi)(\underline{\omega})) = 0$$

holds for any function $\varphi$ depending on $\omega_\ell, \omega_{\ell+1}, \ldots, \omega_\mathfrak{r}$. We begin the computation of the left-hand side by changing variables in the (product-)expectation in order to obtain only $\varphi(\underline{\omega})$:

$$\begin{aligned}
&\mathbf{E}^{(\ell,\mathfrak{r},\theta)}((G^{(\ell,\mathfrak{r},\theta)}\varphi)(\underline{\omega})) \\
&= \mathbf{E}^{(\ell,\mathfrak{r},\theta)}\Bigg\{\sum_{i=\ell}^{\mathfrak{r}-1}[r(\omega_i) + r(-\omega_{i+1})] \cdot [\varphi(\underline{\omega}^{(i,i+1)}) - \varphi(\underline{\omega})] \\
&\qquad + [e^\theta + r(-\omega_\ell)] \cdot [\varphi(\omega_\ell + 1, \ldots) - \varphi(\underline{\omega})] \\
&\qquad + [e^{-\theta} + r(\omega_\mathfrak{r})] \cdot [\varphi(\ldots, \omega_\mathfrak{r} - 1) - \varphi(\underline{\omega})]\Bigg\} \\
&= \mathbf{E}^{(\ell,\mathfrak{r},\theta)}\Bigg\{\Bigg[\sum_{i=\ell}^{\mathfrak{r}-1}\bigg([r(\omega_i + 1) \\
&\qquad + r(-\omega_{i+1} + 1)] \cdot \frac{\mu^{(\theta)}(\omega_i + 1) \cdot \mu^{(\theta)}(\omega_{i+1} - 1)}{\mu^{(\theta)}(\omega_i) \cdot \mu^{(\theta)}(\omega_{i+1})} \\
&\qquad\qquad - r(\omega_i) - r(-\omega_{i+1})\bigg) \\
&\qquad + [e^\theta + r(-\omega_\ell + 1)] \cdot \frac{\mu^{(\theta)}(\omega_\ell - 1)}{\mu^{(\theta)}(\omega_\ell)} - [e^\theta + r(-\omega_\ell)] \\
&\qquad + [e^{-\theta} + r(\omega_\mathfrak{r} + 1)] \cdot \frac{\mu^{(\theta)}(\omega_\mathfrak{r} + 1)}{\mu^{(\theta)}(\omega_\mathfrak{r})} - [e^{-\theta} + r(\omega_\mathfrak{r})]\Bigg] \\
&\qquad\qquad\qquad\qquad\qquad\qquad\qquad\qquad\qquad \times \varphi(\underline{\omega})\Bigg\}.
\end{aligned}$$

We can continue by using the special form 2.5 of $\mu^{(\theta)}$ and then 2.1 of $r$:

$$\begin{aligned}
&\mathbf{E}^{(\ell,\mathfrak{r},\theta)}((G^{(\ell,\mathfrak{r},\theta)}\varphi)(\underline{\omega})) \\
&= \mathbf{E}^{(\ell,\mathfrak{r},\theta)}\Bigg\{\Bigg[\sum_{i=\ell}^{\mathfrak{r}-1}\bigg([r(\omega_i + 1) + r(-\omega_{i+1} + 1)]
\end{aligned}$$

$$\times \frac{r(\omega_{i+1})}{r(\omega_i + 1)} - r(\omega_i) - r(-\omega_{i+1})\Big)$$

$$+ [e^\theta + r(-\omega_\ell + 1)] \cdot \frac{r(\omega_\ell)}{e^\theta} - [e^\theta + r(-\omega_\ell)]$$

$$+ [e^{-\theta} + r(\omega_{\mathfrak{r}} + 1)] \cdot \frac{e^\theta}{r(\omega_{\mathfrak{r}} + 1)} - [e^{-\theta} + r(\omega_{\mathfrak{r}})]\Bigg] \cdot \varphi(\underline{\omega})\Bigg\}$$

$$= \mathbf{E}^{(\ell,\mathfrak{r},\theta)}\Bigg\{\Bigg[\sum_{i=\ell}^{\mathfrak{r}-1}(r(\omega_{i+1}) + r(-\omega_i) - r(\omega_i) - r(-\omega_{i+1}))$$

$$+ r(\omega_\ell) + e^{-\theta} - e^\theta - r(-\omega_\ell) + r(-\omega_{\mathfrak{r}})$$

$$+ e^\theta - e^{-\theta} - r(\omega_{\mathfrak{r}})\Bigg] \cdot \varphi(\underline{\omega})\Bigg\} = 0,$$

which completes the proof. □

3.3. *An explicit construction.* For some arguments it is convenient to have an explicit construction of the finite-volume processes in terms of initial configurations and Poisson clocks. By taking a limit, we can then define the infinite process as a measurable function of these ingredients. Let $\underline{N} = (N_i)_{i \in \mathbb{Z}}$ be an i.i.d. collection of homogeneous rate 1 Poisson point processes on $\mathbb{R}_+^2$. A generic point of $\mathbb{R}_+^2$ is $(t,y)$. The first coordinate $t$ represents time and the second coordinate $y$ represents intensity of rate. As usual, we discard a null set of realizations $\underline{N}$ so that we can assume $N_i(B) < \infty$ for all bounded Borel sets $B \subseteq \mathbb{R}_+^2$ and for all $i$, and

$$\sum_{i\in\mathbb{Z}} N_i(\{t\} \times \mathbb{R}_+) \leq 1 \quad \text{and} \quad \sum_{i\in\mathbb{Z}} N_i(\mathbb{R}_+ \times \{y\}) \leq 1 \qquad \text{for all } (t,y) \in \mathbb{R}_+^2.$$

Fix $-\infty < \ell < \mathfrak{r} < \infty$, and consider initial states $\underline{\omega} \in \mathbb{Z}^\mathbb{Z}$ and $\underline{h} = \{h_i\}_{i\in\mathbb{Z}}$ such that (2.4) holds. The growth rate for column $i$ as a function of the increments $\underline{\omega}$ is given by

$$(3.3) \qquad r_i(\underline{\omega}) := \begin{cases} r(\omega_i) + r(-\omega_{i+1}), & \text{for } \ell \leq i < \mathfrak{r}, \\ 0, & \text{otherwise,} \end{cases}$$

for the $[\ell, \mathfrak{r}]$-monotone process, and

$$(3.4) \qquad r_i(\underline{\omega}) := \begin{cases} e^{\theta_1} + r(-\omega_\ell), & \text{for } i = \ell - 1, \\ r(\omega_i) + r(-\omega_{i+1}), & \text{for } \ell \leq i < \mathfrak{r}, \\ e^{-\theta_1} + r(\omega_{\mathfrak{r}}), & \text{for } i = \mathfrak{r}, \\ 0, & \text{otherwise,} \end{cases}$$

for the $(\ell, \mathfrak{r}, \theta)$ process.

The processes are now constructed by the standard Markov chain scheme that proceeds from jump to jump. We illustrate briefly for the process $\underline{\omega}^{[\ell,\mathfrak{r}]}(t)$. Set $\tau_0 = 0$, $\underline{\omega}_{(0)} = \underline{\omega}$ and

$$\tau_1 = \inf\left\{t > 0 : \sum_{i=\ell}^{\mathfrak{r}-1} N_i((0,t] \times (0, r_i(\underline{\omega}))) = 1\right\}.$$

Define the first segment of the process by $\underline{\omega}^{[\ell,\mathfrak{r}]}(t) = \underline{\omega}_{(0)}$ for $\tau_0 \le t < \tau_1$, and then $\underline{\omega}_{(1)} = \underline{\omega}_{(0)}^{(j_1, j_1+1)}$, where $j_1$ is the index for which $N_{j_1}((0,\tau_1] \times (0, r_{j_1}(\underline{\omega}))) = 1$. Continue inductively, with

$$\tau_k = \inf\left\{t > \tau_{k-1} : \sum_{i=\ell}^{\mathfrak{r}-1} N_i((\tau_{k-1}, \tau_k] \times (0, r_i(\underline{\omega}_{(k-1)}))) = 1\right\},$$

$\underline{\omega}^{[\ell,\mathfrak{r}]}(t) = \underline{\omega}_{(k-1)}$ for $\tau_{k-1} \le t < \tau_k$, $\underline{\omega}_{(k)} = \underline{\omega}_{(k-1)}^{(j_k, j_k+1)}$, and $j_k$ the index for which $N_{j_k}((\tau_{k-1}, \tau_k] \times (0, r_{j_k}(\underline{\omega}_{(k-1)}))) = 1$. We already showed that these finite-volume processes are well defined for all time, so no blow-up happens.

For the future we observe a regularity property of the construction. For this purpose, let us discard another null set of realizations of $\underline{N}$ so we can assume

(3.5) $$\sum_{i\in\mathbb{Z}} \sum_{k,l\in\mathbb{Z}} N_i(\mathbb{R}_+ \times \{r(k) + r(l)\}) = 0.$$

With this assumption, $\underline{\omega}^{[\ell,\mathfrak{r}]}(\cdot)$ is a continuous $D_\Omega$-valued function of $(\underline{\omega}, \underline{N})$. For this claim, put the product topology on $\Omega$ and the product of vague topologies on the point measure configurations $\underline{N}$. $D_\Omega$ denotes the standard space of $\Omega$-valued cadlag paths. To see this continuity on a fixed time interval $[0,T]$, fix $(\underline{\omega}, \underline{N})$ and suppose $(\underline{\omega}', \underline{N}')$ is close enough to $(\underline{\omega}, \underline{N})$. Then $\omega_i' = \omega_i$ for $\ell \le i < \mathfrak{r}$, and in a large enough bounded set of $\mathbb{R}_+^2$ points of $N_i'$ and $N_i$ ($\ell \le i < \mathfrak{r}$) correspond to each other bijectively and are only a small $\varepsilon > 0$ apart. By considering each jump in turn, we see that the same sequence of states $\underline{\omega}_{(0)}, \underline{\omega}_{(1)}, \underline{\omega}_{(2)}, \ldots$ is produced, and the jump times differ by at most $\varepsilon$. How small $\varepsilon$ must be taken is determined by how close together the original jump times $\tau_k$ are up to time $T$, and how close to the boundaries $r_i(\underline{\omega}_{(k)})$ the Poisson points of $N_i$ appear up to time $T$. This argument works because no Poisson point appears exactly on a boundary. This was the purpose of assumption (3.5).

3.4. *Coupling the processes.* In what follows, we will use a common realization of the Poisson processes for different monotone or $(\ell, \mathfrak{r}, \theta)$ processes. Then the construction, together with attractivity of the processes, gives the following coupling. Fix integers $\mathcal{L} \le \ell < \mathfrak{r} \le \mathcal{R}$ and any initial configurations $\underline{\zeta}(0), \underline{\omega}(0) \in \Omega$. Process $\underline{\zeta}(\cdot)$ evolves either according to generator $L^{[\mathcal{L},\mathcal{R}]}$ (2.7) or $G^{(\mathcal{L},\mathcal{R},\theta)}$ (3.1). The other one $\underline{\omega}(\cdot)$ evolves according to $L^{[\ell,\mathfrak{r}]}$. The coupling that results from using common Poisson processes in the construction is summarized in rate tables, rather than by writing explicitly the complicated generator for the coupled process. In these tables, we assume generator $L^{[\mathcal{L},\mathcal{R}]}$ for $\underline{\zeta}$, rather than $G^{(\mathcal{L},\mathcal{R},\theta)}$; the modification needed for $G^{(\mathcal{L},\mathcal{R},\theta)}$ is indicated before Lemma 3.2.

TABLE 1
*Rate for bricklayers at sites $\mathcal{L} \leq i < \ell$ or $\mathfrak{r} \leq i < \mathcal{R}$ to lay brick to their right*

| With rate | $h_i$ | $g_i$ | $d_i$ | $d_{i+1}$ |
|---|---|---|---|---|
| $r(\zeta_i)$ | | ↑ | ↓ | ↑ |

As remarked in Section 2.1, it is enough to describe the rates for individual bricklayers standing at the sites. Given configurations $\underline{\zeta}$ and $\underline{\omega}$, the move of laying a brick to the left is independent of laying one to the right for each bricklayer, hence, we give separate tables for these steps. We also distinguish between sites where both $\underline{\zeta}$ and $\underline{\omega}$ grow, and sites where only columns of $\underline{\zeta}$ grow. We denote the height of the column between $i$ and $i+1$ by $g_i(t)$ [or $h_i(t)$] for the $\underline{\zeta}$ process (or the $\underline{\omega}$ process, resp.). We introduce the notation $d_i := \zeta_i - \omega_i$. As we are essentially interested in this quantity, we give separate columns in the tables as well to describe its behavior. Finally, let ↑ (↓) mean that the corresponding quantity increases (decreases, resp.) by one. For example, $h_i \uparrow$ represents the jump $\underline{\omega} \to \underline{\omega}^{(i,i+1)}$.

Marginally, the processes $\underline{\zeta}(\cdot)$ and $\underline{\omega}(\cdot)$ evolve according to generators $L^{[\mathcal{L},\mathcal{R}]}$ and $L^{[\ell,\mathfrak{r}]}$, respectively. We say that $d_i$ *second class particles are present at site $i$* if $d_i > 0$, and $-d_i$ *second class antiparticles are present at site $i$* if $d_i < 0$.

In case $\underline{\zeta}$ evolves according to $G^{(\mathcal{L},\mathcal{R},\theta)}$, then Table 1 is also valid for $i = \mathcal{R}$, but with rate $r(\zeta_{\mathcal{R}}) + e^{-\theta}$, and Table 2 is valid for $i = \mathcal{L}$ with rate $r(-\zeta_{\mathcal{L}}) + e^{\theta}$.

LEMMA 3.2. *Fix $\mathcal{L} \leq \ell < \mathfrak{r} \leq \mathcal{R}$ and couple an $[\mathcal{L},\mathcal{R}]$-monotone process or $(\mathcal{L},\mathcal{R},\theta)$-process $\underline{\zeta}$ with an $[\ell,\mathfrak{r}]$-monotone process $\underline{\omega}$. Assume that initially $h_i(0) \leq g_i(0)$ for each $i \in \mathbb{Z}$. Then $h_i(t) \leq g_i(t)$ for all $i \in \mathbb{Z}$ and $t \geq 0$.*

TABLE 2
*Rate for bricklayers at sites $\mathcal{L} < i \leq \ell$ or $\mathfrak{r} < i \leq \mathcal{R}$ to lay brick to their left*

| With rate | $h_{i-1}$ | $g_{i-1}$ | $d_{i-1}$ | $d_i$ |
|---|---|---|---|---|
| $r(-\zeta_i)$ | | ↑ | ↓ | ↑ |

TABLE 3
*Rates for bricklayers at sites $\ell \leq i < \mathfrak{r}$ to lay brick to their right*

| With rate | $h_i$ | $g_i$ | $d_i$ | $d_{i+1}$ |
|---|---|---|---|---|
| $[r(\zeta_i) - r(\omega_i)]^+$ | | ↑ | ↓ | ↑ |
| $[r(\omega_i) - r(\zeta_i)]^+$ | ↑ | | ↑ | ↓ |
| $\min[r(\zeta_i), r(\omega_i)]$ | ↑ | ↑ | | |

PROOF. Let us see how the condition stated could be first violated. The step which first violates $h_i \leq g_i$ must be a growth of $h_i$, in a situation when $h_i = g_i$ holds. This step implies the growth of $d_i$ and the decrease of $d_{i+1}$, that is, an $(\uparrow, \downarrow)$ change in $(d_i, d_{i+1})$. In the previous tables, only the second line of Table 3 and the first line of Table 4 contain such a move. However, we are after the first case when the condition would be violated, thus by

$$h_{i-1} \leq g_{i-1}, \qquad h_i = g_i, \qquad h_{i+1} \leq g_{i+1},$$

we have $\zeta_i \geq \omega_i$ and $\zeta_{i+1} \leq \omega_{i+1}$. Under this condition, both moves mentioned above have zero rates, hence, the condition can a.s. never be violated. □

LEMMA 3.3. *Fix $\mathcal{L} < \ell$ and $\mathfrak{r}$. Let $h_i(t)$ be the heights of a process $\underline{\omega}(t)$ evolving according to the $[\ell, \mathfrak{r}]$-monotone evolution, and $g_i(t)$ the heights of $\underline{\zeta}(t)$ evolving according to the $[\mathcal{L}, \mathfrak{r}]$-monotone evolution, coupled as described above. Assume also $h_i(0) = g_i(0)$ for each site $i \in \mathbb{Z}$. Then for all later times, in addition to Lemma 3.2, we also have*

$$\omega_i(t) \leq \zeta_i(t) \tag{3.6}$$

*for each $\ell \leq i \leq \mathfrak{r}$.*

PROOF. Clearly, the statement holds initially. Steps which could violate it are the ones containing $d_i \downarrow$ in the coupling tables. These steps do not apply to sites $\ell \leq i \leq \mathfrak{r}$ in Tables 1 and 2. While (3.6) holds, such steps only apply in Tables 3 and 4 when $\omega_i < \zeta_i$. However, in this case, (3.6) will not be violated. □

TABLE 4
*Rates for bricklayers at sites $\ell < i \leq \mathfrak{r}$ to lay brick to their left*

| With rate | $h_{i-1}$ | $g_{i-1}$ | $d_{i-1}$ | $d_i$ |
|---|---|---|---|---|
| $[r(-\omega_i) - r(-\zeta_i)]^+$ | ↑ | | ↑ | ↓ |
| $[r(-\zeta_i) - r(-\omega_i)]^+$ | | ↑ | ↓ | ↑ |
| $\min[r(-\omega_i), r(-\zeta_i)]$ | ↑ | ↑ | | |

LEMMA 3.4. *Fix $\ell$ and $\mathfrak{r} < \mathcal{R}$. Let $h_i(t)$ be the heights of a process $\underline{\omega}(t)$ evolving according to the $[\ell, \mathfrak{r}]$-monotone evolution, and $g_i(t)$ the heights of $\underline{\zeta}(t)$ evolving according to the $[\ell, \mathcal{R}]$-monotone evolution, coupled as described above. Assume also $h_i(0) = g_i(0)$ for each site $i \in \mathbb{Z}$. Then for all later times, in addition to Lemma 3.2, we also have*

$$\omega_i(t) \geq \zeta_i(t) \tag{3.7}$$

*for each $\ell \leq i \leq \mathfrak{r}$.*

PROOF. Clearly, the statement holds initially. Steps which could violate it are the ones containing $d_i \uparrow$ in the coupling tables. These steps do not apply to sites $\ell \leq i \leq \mathfrak{r}$ in Tables 1 and 2. While (3.7) holds, such steps only apply in Tables 3 and 4 when $\omega_i > \zeta_i$. However, in this case, (3.7) will not be violated. □

**4. The infinite volume limit.** In this section we make use of the connection between the monotone process and the $(\ell, \mathfrak{r}, \theta)$-process. While we have monotonicity described in Lemma 3.2 for the monotone process, we have stationarity for the $(\ell, \mathfrak{r}, \theta)$-process, which gives us stochastic bounds for its growth. Since these bounds are independent of the size $\mathfrak{r} - \ell$ of the $(\ell, \mathfrak{r}, \theta)$-process, our goal is now to construct and then to dominate the limit of the monotone processes by these bounds.

PROOF OF THEOREM 2.2. For each $\ell$ and $\mathfrak{r}$, we use a common realization of the Poisson processes in the construction of Section 3.3. Then Lemma 3.2 shows that $h_i^{[\ell, \mathfrak{r}]}(t)$ is monotone increasing with the extension of the interval $[\ell, \mathfrak{r}]$. Hence, the limits exist a.s. for any initial configuration. □

4.1. *Starting from good distributions.* Recall the definition of a good measure from Section 2.3.

LEMMA 4.1. *Let $\underline{\zeta}(0)$ be distributed according to a good measure $\underline{\pi}$ with parameters $\theta_1 < \theta_2$, and let it evolve according to the $(\ell, \mathfrak{r}, \theta_1)$-evolution. Then there are processes $\underline{\eta}(\cdot)$ and $\underline{\xi}(\cdot)$ with the following properties: $\eta_i(t) \leq \zeta_i(t) \leq \xi_i(t)$ for each site $\ell \leq i \leq \mathfrak{r}$. For each time $t$, the variables $\{\eta_i(t)\}_{\ell \leq i \leq \mathfrak{r}}$ are i.i.d. with common marginal $\mu^{(\theta_1)}$, and the variables $\{\xi_i(t)\}_{\ell \leq i \leq \mathfrak{r}}$ are i.i.d. with common marginal $\mu^{(\theta_2)}$.*

PROOF. By the coupling assumed in the definition of $\underline{\pi}$, we couple $\underline{\zeta}$ with the $(\ell, \mathfrak{r}, \theta_2)$-process $\underline{\xi}$ started from initial distribution $\underline{\mu}^{(\theta_2)}$, and the $(\ell, \mathfrak{r}, \theta_1)$-process $\underline{\eta}$ started from $\underline{\mu}^{(\theta_1)}$, such that $\xi_i(0) \geq \zeta_i(0) \geq \eta_i(0)$ holds initially for each site $\ell \leq i \leq \mathfrak{r}$. Then the basic coupling between these processes, in a very similar way as in Section 3.4, takes care of the inequalities. We show the detailed coupling tables Appendix B. □

LEMMA 4.2. *Let $\underline{\zeta}(0)$ be distributed according to a good measure $\underline{\pi}$ with parameters $\theta_1 < \theta_2$, and let it evolve according to the $(\ell, \mathfrak{r}, \theta_1)$-evolution. Then the columns satisfy*

$$\mathbf{E}[g_i(t) - g_i(0)] \leq t \cdot (e^{\theta_2} + e^{-\theta_1}) \qquad \text{for all } \ell - 1 \leq i \leq \mathfrak{r} \text{ and } 0 \leq t < \infty.$$

PROOF. Consider the growth rate $r_i(\underline{\zeta})$ for column $i$, defined as (3.4) in the $(\ell, \mathfrak{r}, \theta_1)$-process $\underline{\zeta}$. Then

$$(4.1) \qquad M_i(t) := g_i(t) - g_i(0) - \int_0^t r_i(\underline{\zeta}(s)) \, ds$$

is a martingale w.r.t. the filtration generated by $(\underline{\zeta}(s))_{0 \leq s \leq t}$ with $M_i(0) = 0$. By the previous lemma, $\zeta_i(t)$ is bounded from below and from above by $\eta_i(t)$ and $\xi_i(t)$, respectively. Due to monotonicity or $r$, this means

$$(4.2) \qquad r_i(\underline{\zeta}) \leq R_i(\underline{\eta}, \underline{\xi}) := \begin{cases} e^{\theta_1} + r(-\eta_\ell), & \text{for } i = \ell - 1, \\ r(\xi_i) + r(-\eta_{i+1}), & \text{for } \ell \leq i < \mathfrak{r}, \\ e^{-\theta_1} + r(\xi_\mathfrak{r}), & \text{for } i = \mathfrak{r}. \end{cases}$$

Hence, by (4.1),

$$M_i(t) \geq g_i(t) - g_i(0) - \int_0^t R_i(\underline{\eta}(s), \underline{\xi}(s)) \, ds.$$

As $\eta_i(t)$ and $\xi_i(t)$ has (stationary) distributions $\mu^{(\theta_1)}$ and $\mu^{(\theta_2)}$, respectively, taking expectation value of this inequality leads to

$$0 \geq \mathbf{E}[g_i(t) - g_i(0)] - t \cdot (e^{\theta_2} + e^{-\theta_1})$$

for columns $\ell - 1 \leq i \leq \mathfrak{r}$, where we used

$$\mathbf{E}^{(\theta_2)}(r(\xi_i)) = e^{\theta_2}, \qquad \mathbf{E}^{(\theta_1)}(r(-\eta_i)) = e^{-\theta_1}. \qquad \square$$

For later use, we prove here a similar bound for the second moment of $g$. More details about this quantity in equilibrium appear in [3].

LEMMA 4.3. *Let $\underline{\zeta}(0)$ be distributed according to a good measure $\underline{\pi}$ with parameter $\theta_1 < \theta_2$, and let it evolve according to the $(\ell, \mathfrak{r}, \theta_1)$-evolution. Then we have finite constants $A$, $B$, $C$ independent of $\ell$ and $\mathfrak{r}$ such that*

$$\mathbf{E}([g_i(t) - g_i(0)]^2) \leq At^2 + Bt + C \qquad \text{for all } \ell - 1 \leq i \leq \mathfrak{r} \text{ and } 0 \leq t < \infty.$$

PROOF. We introduce

$$\widetilde{g}_i(t) := g_i(t) - g_i(0).$$

With the notation of (3.4), the quantity

$$N_i(t) := \widetilde{g}_i^2(t) - \int_0^t [2\widetilde{g}_i(s) \cdot r_i(s) + r_i(s)] \, ds$$

is a martingale. Hence,

$$\mathbf{E}[\widetilde{g}_i^2(t)] = \int_0^t (2\mathbf{E}[\widetilde{g}_i(s) \cdot r_i(s)] + \mathbf{E}[r_i(s)]) \, ds$$
$$\leq \int_0^t (2\sqrt{\mathbf{E}[\widetilde{g}_i^2(s)]} \cdot \sqrt{\mathbf{E}[r_i^2(s)]} + \mathbf{E}[r_i(s)]) \, ds.$$

Using the bounds given in (4.2) and stationarity of $\underline{\eta}$ and $\underline{\xi}$, the second moments of the rates can be dominated, and we can write

$$\mathbf{E}[\widetilde{g}_i^2(t)] \leq 2\sqrt{\mathbf{E}[R_i^2]} \int_0^t \sqrt{\mathbf{E}[\widetilde{g}_i^2(s)]} \, ds + t \cdot \mathbf{E}[R_i].$$

The quantity $\mathbf{E}[R_i^2]$ contains second moments of the rates w.r.t. $\mu^{(\theta_1)}$ or $\mu^{(\theta_2)}$ (2.5). The existence of these moments can easily be shown, but we do not obtain an explicit formula for them. The left-hand side of the previous equation can be bounded from above by a solution of

$$\frac{d}{dt}x(t) = 2\sqrt{\mathbf{E}[R_i^2]} \cdot \sqrt{x(t)} + \mathbf{E}[R_i].$$

Whenever the increasing function $x(t)$ reaches value 1, its further evolution is dominated by a solution of

$$\frac{d}{dt}y(t) = 3\sqrt{\mathbf{E}[R_i^2]} \cdot \sqrt{y(t)}, \qquad \text{that is, } y(t) = \frac{9}{4}\mathbf{E}[R_i^2] \cdot (t - t_0)^2$$

for some $t_0$, which proves the statement. □

PROOF OF THEOREM 2.3. By Lemma 3.2, for each starting configuration of the $[\ell, \mathfrak{r}]$-monotone process $\underline{\omega}(t)$, there is a bounding $(\ell, \mathfrak{r}, \theta_1)$-process $\underline{\zeta}(t)$ with column heights $g^{[\ell, \mathfrak{r}]}$, for which $\underline{\zeta}(0) = \underline{\omega}(0)$ and

$$h_i^{[\ell,\mathfrak{r}]}(0) = g_i^{[\ell,\mathfrak{r}]}(0), \qquad h_i^{[\ell,\mathfrak{r}]}(t) \leq g_i^{[\ell,\mathfrak{r}]}(t)$$

holds ($\ell \leq i < \mathfrak{r}$). If $\underline{\omega}(0)$ and thus $\underline{\zeta}(0)$ is distributed according to $\underline{\pi}$, then Lemma 4.2 leads to the inequality

(4.3) $\quad \mathbf{E}[h_i^{[\ell,\mathfrak{r}]}(t) - h_i^{[\ell,\mathfrak{r}]}(0)] \leq \mathbf{E}[g_i^{[\ell,\mathfrak{r}]}(t) - g_i^{[\ell,\mathfrak{r}]}(0)] \leq t \cdot (e^{\theta_2} + e^{-\theta_1}).$

Taking liminf of (4.3) leads to similar bounds for $h_i(t)$ via Fatou's lemma. □

4.2. *Conditional coupling: starting from a fixed state.* So far we showed stochastic bounds for the process started from good distributions. Here we apply these bounds to processes started from deterministic initial states $\underline{\omega} \in \widetilde{\Omega}$ of (2.9).

Throughout the rest of the paper, we will use the following conditional coupling construction. Fix an initial state $\underline{\omega} \in \widetilde{\Omega}$. Then by Lemma A.1, there are parameters $\theta_2 > \theta_1$ such that

(4.4) $$\mathbf{E}^{(\theta_1)}(z) < -K_{\underline{\omega}} < K_{\underline{\omega}} < \mathbf{E}^{(\theta_2)}(z).$$

Let $\underline{\pi}$ be the product of marginals

(4.5) $$\pi_i(z) = \begin{cases} \mu^{(\theta_2)}(z), & \text{for } i \leq 0, \\ \mu^{(\theta_1)}(z), & \text{for } i \geq 1. \end{cases}$$

Let $\underline{\zeta}$ have initial distribution $\underline{\pi}$, with initial height $g_0(0) = 0$ at the origin. We define the set

(4.6) $$\mathcal{A}^{\underline{\omega}} := \{\widehat{\underline{\zeta}} \in \Omega : \widehat{g}_i \geq h_i \text{ for all } i \in \mathbb{Z}\}$$

with column heights $h_i$ of $\underline{\omega}$ and $\widehat{g}_i$ of $\widehat{\underline{\zeta}}$. Here (2.4) describes the connection between the $h_i$'s and $\underline{\omega}$, and a similar relation is assumed between the $\widehat{g}_i$'s and $\widehat{\underline{\zeta}}$. With $\underline{\omega}$ fixed, the event $\underline{\zeta} \in \mathcal{A}^{\underline{\omega}}$ might or might not happen when distributing the $\underline{\zeta}$ process initially according to $\underline{\pi}$.

LEMMA 4.4.
$$\underline{\pi}\{\underline{\zeta} \in \mathcal{A}^{\underline{\omega}}\} > 0$$
for each $\underline{\omega} \in \widetilde{\Omega}$.

PROOF. By the assumption on $\underline{\omega}$, there is a bound $K$ and a number $N > 0$ such that

$$K > \frac{1}{|i|} \sum_{j=i+1}^{0} \omega_j = \frac{h_i}{|i|} \qquad \text{for } i < -N,$$

$$K > \frac{1}{i} \sum_{j=1}^{i} (-\omega_j) = \frac{h_i}{i} \qquad \text{for } i > N.$$

It is clear that, with positive probability, $g_i \geq h_i$ happens for all $-N \leq i \leq N$. We show that, with positive probability, this also happens simultaneously for $i < -N$ and $i > N$. For positive $i$'s, due to the previous inequalities, it is enough to show that, with positive probability, $g_i - K \cdot i \geq 0$ for all $i > N$. But the latter is a random walk with drift, of which the increments $-\zeta_i - K$ are i.i.d. random variables and have expectation $-\mathbf{E}^{(\theta_1)}(z) - K > 0$ by (4.4). This random walk has positive probability of never returning to zero. A similar argument works for $i < -N$. □

We can now consider the measure $(\underline{\pi} | \mathcal{A}^{\underline{\omega}})$, which is the measure $\underline{\pi}$ conditioned on the event $\underline{\zeta} \in \mathcal{A}^{\underline{\omega}}$. Probabilities and expectations w.r.t. this conditional measure will be denoted by $\mathbf{P}^{(\underline{\pi} | \mathcal{A}^{\underline{\omega}})}$ and $\mathbf{E}^{(\underline{\pi} | \mathcal{A}^{\underline{\omega}})}$, respectively.

Let $\underline{\zeta}^{(\ell,\mathfrak{r},\theta_1)}$ start from $\underline{\zeta}$ in distribution $\underline{\pi}$, and let it evolve according to the $(\ell,\mathfrak{r},\theta_1)$-evolution, while $\underline{\omega}^{[\ell,\mathfrak{r}]}$, started from the fixed state $\underline{\omega}$, evolves according to the $[\ell,\mathfrak{r}]$-monotone evolution. We couple these two processes as described in Section 3.4. Then:

- the marginal evolution of $\underline{\omega}^{[\ell,\mathfrak{r}]}$ is not influenced by the occurrence of the event $\underline{\zeta} \in \mathcal{A}^{\underline{\omega}}$, since the initial state $\underline{\omega}$ was fixed,
- if $\underline{\zeta} \in \mathcal{A}^{\underline{\omega}}$ happened initially, then, by Lemma 3.2, $h_i^{[\ell,\mathfrak{r}]}(t) \leq g_i^{(\ell,\mathfrak{r},\theta_1)}(t)$ holds for each site $i$ and any time $t > 0$.

We call this construction a *conditional coupling*.

PROOF OF THEOREM 2.4.   By the conditional coupling described above, we have the following:

$$\mathbf{E}^{\underline{\omega}}(h_i^{[\ell,\mathfrak{r}]}(t)) = \mathbf{E}^{(\underline{\pi}|\mathcal{A}^{\underline{\omega}})}(h_i^{[\ell,\mathfrak{r}]}(t)) \leq \mathbf{E}^{(\underline{\pi}|\mathcal{A}^{\underline{\omega}})}(g_i^{[\ell,\mathfrak{r}]}(t))$$
$$= \mathbf{E}^{(\underline{\pi}|\mathcal{A}^{\underline{\omega}})}(g_i^{[\ell,\mathfrak{r}]}(t) - g_i(0)) + \mathbf{E}^{(\underline{\pi}|\mathcal{A}^{\underline{\omega}})}(g_i(0))$$
$$= \frac{\mathbf{E}^{\underline{\omega}}([g_i^{[\ell,\mathfrak{r}]}(t) - g_i(0)] \cdot \mathbf{1}\{\underline{\zeta} \in \mathcal{A}^{\underline{\omega}}\})}{\underline{\pi}\{\underline{\zeta} \in \mathcal{A}^{\underline{\omega}}\}} + \frac{\mathbf{E}^{\underline{\omega}}([g_i(0)] \cdot \mathbf{1}\{\underline{\zeta} \in \mathcal{A}^{\underline{\omega}}\})}{\underline{\pi}\{\underline{\zeta} \in \mathcal{A}^{\underline{\omega}}\}}$$
$$\leq \frac{\mathbf{E}([g_i^{[\ell,\mathfrak{r}]}(t) - g_i(0)])}{\underline{\pi}\{\underline{\zeta} \in \mathcal{A}^{\underline{\omega}}\}} + \frac{\mathbf{E}(|g_i(0)|)}{\underline{\pi}\{\underline{\zeta} \in \mathcal{A}^{\underline{\omega}}\}}.$$

For the first term, we can apply Lemma 4.2 and finally write

$$\mathbf{E}^{\underline{\omega}}(h_i^{[\ell,\mathfrak{r}]}(t)) \leq t \cdot \frac{e^{\theta_2} + e^{-\theta_1}}{\underline{\pi}\{\underline{\zeta} \in \mathcal{A}^{\underline{\omega}}\}} + \frac{\mathbf{E}(|g_i(0)|)}{\underline{\pi}\{\underline{\zeta} \in \mathcal{A}^{\underline{\omega}}\}}.$$

Note that this bound does not depend on $\ell$ and $\mathfrak{r}$ as long as $\ell < i < \mathfrak{r}$. As in Theorem 2.3, we can finish the proof by monotonicity of $h_i^{[\ell,\mathfrak{r}]}$ in $\ell$, $\mathfrak{r}$ and by Fatou's lemma. In the second term, $|g_i(0)|$ is bounded from above by the sum of $|i|$ many $|\zeta(0)|$'s, initially in i.i.d. distribution, hence, the expectation is linear in $|i|$.  □

4.3. *A Markov process of heights with cadlag paths.*  We finish this section by establishing that when the initial increment configuration $\underline{\omega}$ lies in the space $\widetilde{\Omega}$ of states with Cesàro-bounded slopes, then the infinite-volume height process $\underline{h}(\cdot)$ is a Markov process with cadlag paths. We extend these properties to $\underline{\omega}(\cdot)$ after we have shown that the space $\widetilde{\Omega}$ is actually closed under the $\underline{\omega}$-evolution.

PROPOSITION 4.5.    *Let $a < b$ in $\mathbb{Z}$, $0 < T < \infty$ a fixed time, and $\underline{h}$ an initial height configuration whose increments satisfy $\underline{\omega} \in \widetilde{\Omega}$. Then there is an*

a.s. finite positive random variable $R = R^\omega(a, b, T)$ such that

$$h_i^{[\ell, \mathfrak{r}]}(t) = h_i(t) \qquad \text{for } a \leq i < b \text{ and } 0 \leq t \leq T,$$

if $\ell < -R < R < \mathfrak{r}$. In particular, the $\mathbb{Z}^{(b-a)}$-valued path $t \mapsto \{h_i(t)\}_{a \leq i < b}$ of the limiting process is a.s. a cadlag function on $[0, T]$.

PROOF. For all $i$ and $\ell < \mathfrak{r}$, $h_i^{[\ell, \mathfrak{r}]}(t) \leq h_i(t)$ a.s. and these are both non-decreasing in $t$. Therefore, for $t \in [0, T]$,

$$\omega_i^{[\ell, \mathfrak{r}]}(t) = h_{i-1}^{[\ell, \mathfrak{r}]}(t) - h_i^{[\ell, \mathfrak{r}]}(t) \leq h_{i-1}(t) - h_i^{[\ell, \mathfrak{r}]}(0) \leq h_{i-1}(T) - h_i(0) \qquad \text{a.s.},$$

and the last bound is a.s. finite. A similar estimate works for $-\omega_i^{[\ell, \mathfrak{r}]}(t)$. Thus,

$$\sup\{|\omega_i^{[\ell, \mathfrak{r}]}(t)| : a \leq i \leq b, t \in [0, T], -\infty < \ell < \mathfrak{r} < \infty\} < \infty \qquad \text{almost surely.}$$

Consequently, the rates $r_i(\underline{\omega}^{[\ell, \mathfrak{r}]}(t))$ of (3.3) are a.s. bounded, uniformly over the relevant range of indices and time. In the coupling construction of Section 3.4 we used common Poisson processes $\underline{N}$ for the $[\ell, \mathfrak{r}]$-monotone processes. The uniform bound on the rates then implies that the family of height processes $\{h_i^{[\ell, \mathfrak{r}]}(t) : a \leq i < b, t \in [0, T], -\infty < \ell < \mathfrak{r} < \infty\}$ is constructed with a single realization of the Poisson processes $\{N_i : a \leq i < b\}$ restricted to a bounded (but random) rectangle $[0, T] \times [0, A] \subseteq \mathbb{R}_+^2$. So this entire height family has a fixed finite (but random) set of possible jump times. Since the limiting heights are now finite as proved above, for each fixed $(t, i)$, the convergence implies that $h_i^{[\ell, \mathfrak{r}]}(t) = h_i(t)$ for large enough $\ell, \mathfrak{r}$. If we take $\ell, \mathfrak{r}$ large enough to have the equality on the set of possible jump times and at the endpoint $T$, then the equality must hold throughout the time interval $[0, T]$. And finally, since each path $t \mapsto h_i^{[\ell, \mathfrak{r}]}(t)$ is a nondecreasing integer-valued cadlag function, this conclusion follows for the limiting path. □

PROPOSITION 4.6. *Let $\underline{h}$ be an initial height configuration whose increments satisfy $\underline{\omega} \in \widetilde{\Omega}$. Then the height process $\underline{h}(t)$ is a Markov process.*

PROOF. We can represent the construction of Section 3.3 in terms of a family of measurable mappings $\Phi_t^{[\ell, \mathfrak{r}]}$ defined on initial height configurations and the subset of well-behaved point measures on $\mathbb{R}_+^2$, so that $h^{[\ell, \mathfrak{r}]}(t) = \Phi_t^{[\ell, \mathfrak{r}]}(\underline{h}, \underline{N})$ for all $\underline{h} \in \mathbb{Z}^{\mathbb{Z}}$. The action of this mapping could be characterized by "freeze the columns outside $[\ell, \mathfrak{r}]$ and evolve the columns over $[\ell, \mathfrak{r}]$ as functions of $\{h_i, N_i : \ell \leq i < \mathfrak{r}\}$ as described before." Then as explained, monotonicity gives the existence of limiting measurable mappings

$$(4.7) \qquad \underline{h}(t) = \Phi_t(\underline{h}, \underline{N}) := \lim_{-\ell, \mathfrak{r} \to \infty} \Phi_t^{[\ell, \mathfrak{r}]}(\underline{h}, \underline{N}).$$

The Markov property follows from showing

(4.8) $$\underline{h}(s+t) = \Phi_t(\underline{h}(s), \theta_s \underline{N}) \qquad \text{for } s, t > 0,$$

where $\theta_s$ translates the Poisson points on $\mathbb{R}_+^2$ by $(-s, 0)$ (in other words, moves the time origin to $s$ and "restarts" the Poisson processes). Since the initial increments $\underline{\omega}$ are assumed to lie in $\widetilde{\Omega}$, the height configuration $\underline{h}(s)$ is finite a.s. and, consequently, the values $\Phi_t(\underline{h}(s), \theta_s \underline{N})$ are well defined by the limiting procedure. Fix an index $i$. We will show

(4.9) $$h_i(s+t) = \Phi_{t,i}(\underline{h}(s), \theta_s \underline{N}),$$

where $\Phi_{t,i}$ is the value of the mapping in the $i$th column. Starting from the right-hand side of (4.9), by the definitions,

$$\Phi_{t,i}(\underline{h}(s), \theta_s \underline{N}) = \lim_{-\ell, \mathfrak{r} \to \infty} \Phi_{t,i}^{[\ell, \mathfrak{r}]}(\Phi_s(\underline{h}, \underline{N}), \theta_s \underline{N}).$$

The mapping $\Phi_{t,i}^{[\ell, \mathfrak{r}]}$ reads its input $\Phi_s(\underline{h}, \underline{N})$ only on the sites in $[\ell, \mathfrak{r} - 1]$. Again by the definition of the limiting mapping and the finiteness of the limit, for each $\ell, \mathfrak{r}$, there exist integers $a = a(\ell, \mathfrak{r}) < \ell < \mathfrak{r} < b(\ell, \mathfrak{r}) = b$ such that

$$\Phi_{s,j}(\underline{h}, \underline{N}) = \Phi_{s,j}^{[a,b]}(\underline{h}, \underline{N}) \qquad \text{for } \ell \leq j < \mathfrak{r}.$$

Consequently,

$$\lim_{-\ell, \mathfrak{r} \to \infty} \Phi_{t,i}^{[\ell, \mathfrak{r}]}(\Phi_s(\underline{h}, \underline{N}), \theta_s \underline{N}) = \lim_{-\ell, \mathfrak{r} \to \infty} \Phi_{t,i}^{[\ell, \mathfrak{r}]}(\Phi_s^{[a,b]}(\underline{h}, \underline{N}), \theta_s \underline{N}).$$

Since $[\ell, \mathfrak{r}] \subseteq [a, b]$, monotonicity and the Markov property of the processes $\underline{h}^{[\ell, \mathfrak{r}]}(\cdot)$ and $\underline{h}^{[a,b]}(\cdot)$ give

$$\Phi_{t,i}(\underline{h}(s), \theta_s \underline{N}) \leq \lim_{-\ell, \mathfrak{r} \to \infty} \Phi_{t,i}^{[a,b]}(\Phi_s^{[a,b]}(\underline{h}, \underline{N}), \theta_s \underline{N})$$

$$= \lim_{-\ell, \mathfrak{r} \to \infty} \Phi_{s+t,i}^{[a,b]}(\underline{h}, \underline{N}) = h_i(s+t)$$

and

$$\Phi_{t,i}(\underline{h}(s), \theta_s \underline{N}) \geq \lim_{-\ell, \mathfrak{r} \to \infty} \Phi_{t,i}^{[\ell, \mathfrak{r}]}(\Phi_s^{[\ell, \mathfrak{r}]}(\underline{h}, \underline{N}), \theta_s \underline{N})$$

$$= \lim_{-\ell, \mathfrak{r} \to \infty} \Phi_{s+t,i}^{[\ell, \mathfrak{r}]}(\underline{h}, \underline{N}) = h_i(s+t).$$

This verifies (4.8). $\square$

Utilizing the jointly measurable representation (4.7), we can define a semigroup $S(t)$ of operators acting on bounded measurable functions $f$ by integrating away the Poisson clocks:

$$S(t)f(\underline{h}) = \int f(\Phi_t(\underline{h}, \underline{N})) \mathbb{P}(d\underline{N}).$$

The semigroup property $S(s+t) = S(s)S(t)$ follows from (4.8).

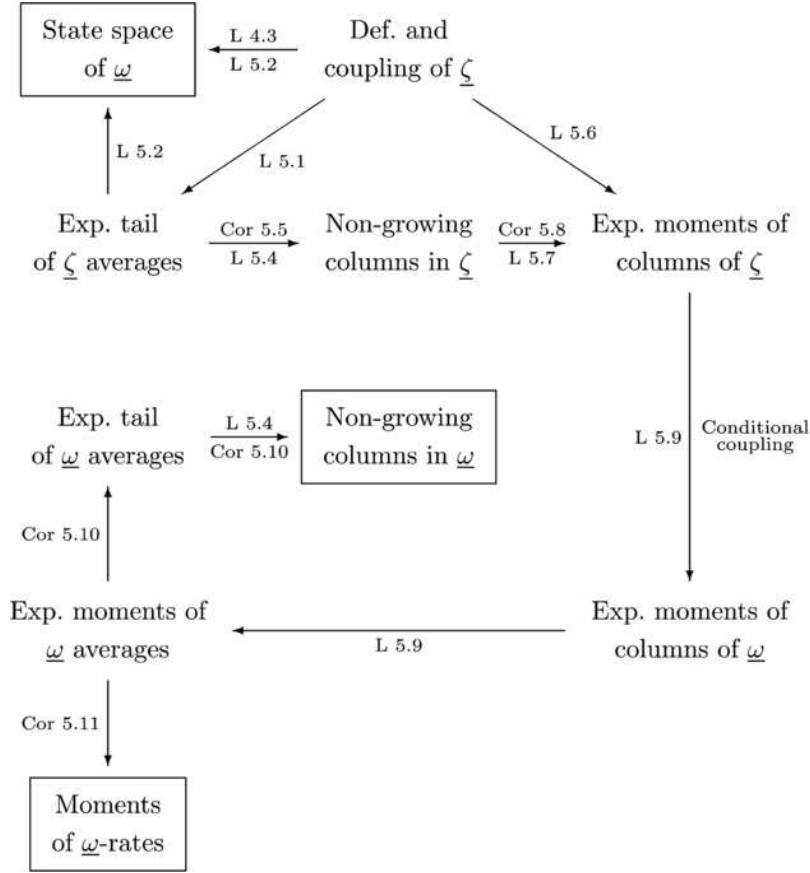

Fig. 2.

**5. Regularity properties.** In the previous section we constructed the process $\underline{\omega}(\cdot)$ as a limit of the monotone processes. In the remaining sections we study this process more closely. The first order of business is to show that the state space $\widetilde{\Omega}$ (2.9) is closed under the dynamics. For later use, we also show that the growth rates of the process are in $L^p$ for any $1 \leq p < \infty$, at least until some time $T = T^{\underline{\omega}}(p)$, and that it is exponentially unlikely for each column in an interval to grow in a small time-range.

The goal is achieved in several steps by first proving these properties for the $(\ell, \mathfrak{r}, \theta)$-process $\underline{\zeta}$, and then using the conditional coupling construction of Section 4.2. The relationships between various statements are illustrated by the diagram below. The boxed statements are the main outcomes of this section and will be used subsequently to study the semigroup and ergodicity (see Figure 2).

5.1. *Markov process of increments.* In this section we complete the basic description of the process $\underline{\omega}(\cdot)$ by proving Theorem 2.5.

LEMMA 5.1. *Let $\underline{\zeta}(t)$ be an $(\ell, \mathfrak{r}, \theta_1)$-process started from distribution $\underline{\pi}$ of 4.5. Then there exists a function $\widetilde{C}(K)$ that depends on $\theta_1$ and $\theta_2$, but not on $t, \ell, \mathfrak{r}$, such that $\widetilde{C}(K) \nearrow \infty$ as $K \nearrow \infty$, and for all $\ell - 1 \leq i < 0 < k \leq \mathfrak{r}$,*

$$\mathbf{P}\left\{\frac{1}{|i|}\sum_{j=i+1}^{0}|\zeta_j(t)| \geq K\right\} \leq e^{-\widetilde{C}(K)|i|} \quad \text{and} \quad \mathbf{P}\left\{\frac{1}{k}\sum_{j=1}^{k}|\zeta_j(t)| \geq K\right\} \leq e^{-\widetilde{C}(K)k}.$$

PROOF. By the coupling in Lemma 4.1, $\eta_j(t) \leq \zeta_j(t) \leq \xi_j(t)$, and the bounds follow because at a fixed time the variables $\{\eta_j(t)\}$ and $\{\xi_j(t)\}$ are i.i.d. with all exponential moments. □

LEMMA 5.2. *Let $\underline{\omega}(t)$ be the limiting process started from $\underline{\omega} \in \widetilde{\Omega}$. Then there is a real $0 < \widehat{K}^{\underline{\omega}} < \infty$ such that*

$$\mathbf{P}^{\underline{\omega}}\left\{\sup_{n \geq 0}\limsup_{i \to \infty}\sup_{t \in [n, n+1]}\frac{1}{i}\sum_{j=-i}^{i}|\omega_j(t)| \leq \widehat{K}^{\underline{\omega}}\right\} = 1.$$

PROOF. Fix an initial state $\underline{\omega} \in \widetilde{\Omega}$ and consider the $[\ell, \mathfrak{r}]$-monotone process $\underline{\omega}^{[\ell, \mathfrak{r}]}$ with $\ell < i < -1$. Let the process $\underline{\zeta}^{(\ell, \mathfrak{r}, \theta_1)}$ be as in the conditional coupling in Section 4.2. Start with the inequality

$$(5.1) \quad \sum_{j=i+1}^{0}|\omega_j^{[\ell,\mathfrak{r}]}(t)| \leq \sum_{j=i+1}^{0}|\omega_j^{[\ell,\mathfrak{r}]}(t) - \zeta_j^{(\ell,\mathfrak{r},\theta_1)}(t)| + \sum_{j=i+1}^{0}|\zeta_j^{(\ell,\mathfrak{r},\theta_1)}(t)|.$$

The first sum on the right is the number of second class particles between $\underline{\zeta}^{(\ell,\mathfrak{r},\theta_1)}$ and $\underline{\omega}^{[\ell,\mathfrak{r}]}$, at sites $i+1, i+2, \ldots, 0$, at time $t$. For any $\ell \leq i < 0 < \mathfrak{r}$, this number can only increase by second class particles jumping over the edge $[i, i+1]$ and $[0,1]$. Any such jump requires at least one of the columns $g^{(\ell,\mathfrak{r},\theta_1)}$ of $\underline{\zeta}^{(\ell,\mathfrak{r},\theta_1)}$ or $h^{[\ell,\mathfrak{r}]}$ of $\underline{\omega}^{[\ell,\mathfrak{r}]}$ over the corresponding edge to grow. Therefore, also taking into consideration $g_0(0) = h_0(0) = 0$ and (2.4),

$$\sum_{j=i+1}^{0}|\omega_j^{[\ell,\mathfrak{r}]}(t) - \zeta_j^{(\ell,\mathfrak{r},\theta_1)}(t)|$$
$$\leq \sum_{j=i+1}^{0}|\omega_j(0) - \zeta_j(0)| + g_i^{(\ell,\mathfrak{r},\theta_1)}(t) - g_i(0) + h_i^{[\ell,\mathfrak{r}]}(t) - h_i(0)$$
$$\quad + g_0^{(\ell,\mathfrak{r},\theta_1)}(t) - g_0(0) + h_0^{[\ell,\mathfrak{r}]}(t) - h_0(0)$$

$$\leq 2 \sum_{j=i+1}^{0} |\omega_j(0)| + 2 \sum_{j=i+1}^{0} |\zeta_j(0)| + g_i^{(\ell,\mathfrak{r},\theta_1)}(t) + h_i^{[\ell,\mathfrak{r}]}(t) + g_0^{(\ell,\mathfrak{r},\theta_1)}(t)$$
$$+ h_0^{[\ell,\mathfrak{r}]}(t).$$

Now consider a range $t_0 \leq t \leq t_1$ of times. By the monotonicity of the heights,

$$\sup_{t_0 \leq t \leq t_1} \sum_{j=i+1}^{0} |\omega_j^{[\ell,\mathfrak{r}]}(t) - \zeta_j^{(\ell,\mathfrak{r},\theta_1)}(t)|$$
$$\leq 2 \sum_{j=i+1}^{0} |\omega_j(0)| + 2 \sum_{j=i+1}^{0} |\zeta_j(0)|$$
$$+ g_i^{(\ell,\mathfrak{r},\theta_1)}(t_1) + h_i^{[\ell,\mathfrak{r}]}(t_1) + g_0^{(\ell,\mathfrak{r},\theta_1)}(t_1) + h_0^{[\ell,\mathfrak{r}]}(t_1).$$

The sum of $\zeta$-increments in (5.1) is handled by introducing

$$v(t) = \sum_{j=i+1}^{0} |\zeta_j^{(\ell,\mathfrak{r},\theta_1)}(t)| - g_i^{(\ell,\mathfrak{r},\theta_1)}(t) - g_0^{(\ell,\mathfrak{r},\theta_1)}(t).$$

Considering the rates of the jumps that increase the quantity, $v(t)$ reveals that the increment $v(t_0 + s) - v(t_0)$ is bounded stochastically by $2Y(s)$, where $Y(\cdot)$ denotes a Poisson process with rate $2(|i| - 1)r(0)$. We have the bound

$$\sup_{t_0 \leq t \leq t_1} \sum_{j=i+1}^{0} |\zeta_j^{(\ell,\mathfrak{r},\theta_1)}(t)|$$
$$\leq \sup_{t_0 \leq t \leq t_1} v(t) + g_i^{(\ell,\mathfrak{r},\theta_1)}(t_1) + g_0^{(\ell,\mathfrak{r},\theta_1)}(t_1)$$
$$\leq \sup_{0 \leq s \leq t_1 - t_0} (v(t_0 + s) - v(t_0)) + \sum_{j=i+1}^{0} |\zeta_j^{(\ell,\mathfrak{r},\theta_1)}(t_0)| - g_i^{(\ell,\mathfrak{r},\theta_1)}(t_0)$$
$$- g_0^{(\ell,\mathfrak{r},\theta_1)}(t_0) + g_i^{(\ell,\mathfrak{r},\theta_1)}(t_1) + g_0^{(\ell,\mathfrak{r},\theta_1)}(t_1)$$
$$\leq \sup_{0 \leq s \leq t_1 - t_0} (v(t_0 + s) - v(t_0)) + 2 \sum_{j=i+1}^{0} |\zeta_j^{(\ell,\mathfrak{r},\theta_1)}(t_0)| - 2g_0^{(\ell,\mathfrak{r},\theta_1)}(t_0)$$
$$+ g_i^{(\ell,\mathfrak{r},\theta_1)}(t_1) + g_0^{(\ell,\mathfrak{r},\theta_1)}(t_1).$$

Combine the estimates to write

$$\sup_{t_0 \leq t \leq t_1} \sum_{j=i+1}^{0} |\omega_j^{[\ell,\mathfrak{r}]}(t)| \leq 2 \sum_{j=i+1}^{0} |\omega_j(0)| + 2 \sum_{j=i+1}^{0} |\zeta_j(0)|$$

$$
\begin{aligned}
(5.2)\quad & + 2 \sum_{j=i+1}^{0} |\zeta_j^{(\ell,\mathfrak{r},\theta_1)}(t_0)| + \sup_{0 \le s \le t_1 - t_0} (v(t_0+s) - v(t_0)) \\
& + h_i^{[\ell,\mathfrak{r}]}(t_1) + 2 g_i^{(\ell,\mathfrak{r},\theta_1)}(t_1) + 2 g_0^{(\ell,\mathfrak{r},\theta_1)}(t_1) \\
& - 2 g_0^{(\ell,\mathfrak{r},\theta_1)}(t_0) + h_0^{[\ell,\mathfrak{r}]}(t_1).
\end{aligned}
$$

We now observe for each term on the right-hand side a bound on the probability of exceeding $K|i|$ that is independent of $\ell, \mathfrak{r}$ and summable in $|i|$, once $K$ is fixed sufficiently large. Then Borel–Cantelli concludes the proof.

By the assumption $\underline{\omega} \in \widetilde{\Omega}$, the first term $\sum |\omega_j(0)|$ is deterministically less than $K|i|$ for a large enough $K$ and $i$, independently of $\ell, \mathfrak{r}$.

To the second and third terms (sums of $\zeta$-variables) apply Lemma 5.1.

As observed above, the fourth term is stochastically dominated by twice a Poisson variable with mean $2(|i|-1)r(0)(t_1 - t_0)$.

For the term $h_i^{[\ell,\mathfrak{r}]}(t_1)$ use conditional coupling. With the event $\underline{\zeta} \in \mathcal{A}^{\underline{\omega}}$ defined in (4.6),

$$
\begin{aligned}
\mathbf{P}^{\underline{\omega}}\{h_i^{[\ell,\mathfrak{r}]}(t_1) > K \cdot |i|\} \\
= \mathbf{P}^{(\underline{\pi}|\mathcal{A}^{\underline{\omega}})}\{h_i^{[\ell,\mathfrak{r}]}(t_1) > K \cdot |i|\} \\
\le \mathbf{P}^{(\underline{\pi}|\mathcal{A}^{\underline{\omega}})}\{g_i^{(\ell,\mathfrak{r},\theta_1)}(t_1) > K \cdot |i|\} \le \frac{\mathbf{P}\{g_i^{(\ell,\mathfrak{r},\theta_1)}(t_1) > K \cdot |i|\}}{\underline{\pi}\{\underline{\zeta} \in \mathcal{A}^{\underline{\omega}}\}} \\
\le \frac{\mathbf{P}\{\sum_{j=i+1}^{0} |\zeta_j^{(\ell,\mathfrak{r},\theta_1)}(t_1)| > K/2 \cdot |i|\}}{\underline{\pi}\{\underline{\zeta} \in \mathcal{A}^{\underline{\omega}}\}} + \frac{\mathbf{P}\{g_0^{(\ell,\mathfrak{r},\theta_1)}(t_1) > (K/2) \cdot |i|\}}{\underline{\pi}\{\underline{\zeta} \in \mathcal{A}^{\underline{\omega}}\}}.
\end{aligned}
$$

As $\underline{\pi}\{\underline{\zeta} \in \mathcal{A}^{\underline{\omega}}\} > 0$ and independent of $\ell$ and $\mathfrak{r}$, this bound is taken care of by the already existing terms in (5.2). Also, the last inequality shows how the term $g_i^{(\ell,\mathfrak{r},\theta_1)}(t_1)$ is reduced to $g_0^{(\ell,\mathfrak{r},\theta_1)}(t_1)$ and another sum of $\zeta$-terms.

The term $g_0^{(\ell,\mathfrak{r},\theta_1)}(t_1)$ of (5.2) is of bounded second moment, uniformly in $\ell$ and $\mathfrak{r}$, by Lemma 4.3. Hence, by Chebyshev's inequality,

$$
\mathbf{P}\{g_0^{(\ell,\mathfrak{r},\theta_1)}(t_1) > K \cdot |i|\} \le \frac{\mathbf{E}(g_0^{(\ell,\mathfrak{r},\theta_1)}(t_1)^2)}{K^2 i^2} \le \frac{A t_1^2 + B t_1 + C}{K^2 i^2},
$$

which is again summable in $i$, uniformly in $\ell$ and $\mathfrak{r}$. The term $g_0^{(\ell,\mathfrak{r},\theta_1)}(t_0)$ is treated similarly.

Finally, for $h_0^{[\ell,\mathfrak{r}]}(t_1)$ repeat the conditional coupling argument to reduce it to $g_0^{(\ell,\mathfrak{r},\theta_1)}(t_1)$.

The estimates above can be summarized as follows: there exist quantities $0 < K < \infty$ and $0 \le a_{|i|} \le 1$ for $i < 0$, independent of $\ell$ and $\mathfrak{r}$, such that

$\sum a_{|i|} < \infty$ and

$$\mathbf{P}^{\underline{\omega}}\left\{\sup_{t_0 \leq t \leq t_1} \frac{1}{|i|} \sum_{j=i+1}^{0} |\omega_j^{[\ell,\mathfrak{r}]}(t)| > 15K\right\} \leq a_{|i|}.$$

Both $K$ and $a_{|i|}$ depend on $\underline{\omega}$, but we leave out the dependence from the notation. By Proposition 4.5, for any fixed $i < 0$ and $t_1 < \infty$,

$$\lim_{-\ell,\mathfrak{r} \to \infty} \mathbf{P}^{\underline{\omega}}\{\ \omega_j(t) = \omega_j^{[\ell,\mathfrak{r}]}(t) \text{ for } i+1 \leq j \leq 0 \text{ and } t \in [0,t_1]\ \} = 1.$$

So we can pass to the limit to obtain the bound for the limiting $\underline{\omega}$-process:

$$\mathbf{P}^{\underline{\omega}}\left\{\sup_{t_0 \leq t \leq t_1} \frac{1}{|i|} \sum_{j=i+1}^{0} |\omega_j(t)| > 15K\right\} \leq a_{|i|}.$$

The first Borel–Cantelli lemma implies that

$$\limsup_{i \to -\infty} \sup_{t_0 \leq t \leq t_1} \frac{1}{|i|} \sum_{j=i+1}^{0} |\omega_j(t)| \leq 15K \qquad \text{almost surely.}$$

Looking back, we see that $K$ depends on time only through the increment $t_1 - t_0$, and this only to bound the fourth term on the right of (5.2). So fix the size of the time increment to $t_1 - t_0 = 1$, and let $t_0$ run over nonnegative integers. A similar argument works for positive $i$-values. $\square$

PROOF OF THEOREM 2.5.  In Section 4.3 we constructed path measures $P^{\underline{h}}$ on $D_{\mathbb{Z}^{\mathbb{Z}}}$ for the height process when the initial configuration $\underline{h}$ has increments $\underline{\omega} \in \widetilde{\Omega}$. Here $D_{\mathbb{Z}^{\mathbb{Z}}}$ is the space of $\mathbb{Z}^{\mathbb{Z}}$-valued cadlag paths. To be somewhat formal about this, we can now define the path measure $\mathbf{P}^{\underline{\omega}}$ of the increment process as the distribution of the process $\underline{\omega}(\cdot) := \theta_{-1}\underline{h}(\cdot) - \underline{h}(\cdot)$ under the measure $P^{\underline{h}^0(\underline{\omega})}$. Here $\theta_1$ shifts spatial index by $(\theta_{-1}h)_i = h_{i-1}$, and $\underline{h}^0(\underline{\omega})$ is the height configuration associated to $\underline{\omega}$ via (2.4) normalized by $h_0 = 0$. Then one can deduce the Markov property and the cadlag property for $\underline{\omega}(\cdot)$ and the measurability of $\underline{\omega} \mapsto \mathbf{P}^{\underline{\omega}}$ from the corresponding properties of $\underline{h}(\cdot)$. Lemma 5.2 completes the proof of Theorem 2.5. $\square$

5.2. *Further regularity properties of the $\underline{\zeta}$ process.*  In this section we derive exponential moment bounds for $\underline{\zeta}(t)$. The key tool is a bound on the probability that a large interval of neighboring columns all grew in a short time. In some lemmas we use the notation $\underline{\omega}(t)$ for both the $(\ell, \mathfrak{r}, \theta)$-process and the $[\ell, \mathfrak{r}]$-monotone process to avoid unnecessary duplication.

The first lemma states an elementary counting argument. We encode an upper bound $z$ on the magnitude of the left and right increment of a column $h_j$ as membership $j \in \mathcal{D}(z)$, as defined below.

LEMMA 5.3.　*Fix real numbers $0 \leq b < 1$ and $0 < \alpha < 1-b$. For integers $i < 0$, $z > 0$, and for $\underline{\omega} \in \mathbb{Z}^{\mathbb{Z}}$, define*

(5.3)　$\mathcal{D}(z) = \mathcal{D}(z, i, b, \underline{\omega}) := \{i+1 \leq j \leq \lfloor bi \rfloor - 1 : |\omega_j| < z \text{ and } |\omega_{j+1}| < z\}.$

*Then there exists an integer $i_0 = i_0(b,\alpha) < 0$ such that the following statement holds for any real $K > 0$. If $i < i_0$ and*

(5.4)　$$\frac{1}{|i|} \sum_{j=i+1}^{0} |\omega_j| < K,$$

*then*

$$|\mathcal{D}(z)| \geq \lfloor \alpha |i| \rfloor \qquad \text{for all } z \geq \frac{3K}{1-b-\alpha}.$$

*The same holds for $i > 0$ with the obvious changes.*

PROOF.　Fix
$$z = \frac{3K}{1-b-\alpha},$$
and define
$$\mathcal{D}_1(z) := \{i+1 \leq j \leq \lfloor bi \rfloor - 1 : |\omega_j| < z\} \quad \text{and}$$
$$\mathcal{D}_2(z) := \{i+1 \leq j \leq \lfloor bi \rfloor - 1 : |\omega_{j+1}| < z\}.$$
We show that

(5.5)　$$|\mathcal{D}_1(z)| \geq \frac{\alpha - b + 1}{2} \cdot |i|$$

whenever $|i| > \frac{12}{1-b-\alpha}$. By contradiction, assuming
$$|\mathcal{D}_1(z)| < \frac{\alpha - b + 1}{2} \cdot |i|$$
implies
$$\lfloor bi \rfloor - 1 + |i| - |\mathcal{D}_1(z)| \geq bi - 2 + |i| - |\mathcal{D}_1(z)| > \frac{1-\alpha-b}{2} \cdot |i| - 2$$
many indices $j$ with $|\omega_j| \geq z$ in the discrete interval $\{i+1 \cdots \lfloor bi \rfloor - 1\}$. The contribution of these to
$$\frac{1}{|i|} \sum_{j=i+1}^{0} |\omega_j|$$
is larger than
$$\frac{1-\alpha-b}{2} \cdot z - \frac{2}{|i|} z = \frac{3}{2} K - \frac{6K}{|i|(1-b-\alpha)} > K$$

whenever $|i| > \frac{12}{1-b-\alpha}$, contradicting (5.4).

Since (5.5) holds also for $\mathcal{D}_2(z)$, we have

$$|\mathcal{D}(z)| = |\mathcal{D}_1(z) \cap \mathcal{D}_2(z)| = |\mathcal{D}_1(z)| + |\mathcal{D}_2(z)| - |\mathcal{D}_1(z) \cup \mathcal{D}_2(z)|$$
$$\geq (\alpha - b + 1)|i| - (\lfloor bi \rfloor - 1 + |i|) > \alpha|i|. \qquad \square$$

The next lemma quantifies this idea: if we have some control of the increments at time $t$, then we also have some control of the growth of columns during an immediately preceding small time interval. The argument is the same for negative and positive indices. We work with negative indices $i$, and define the event

(5.6) $\qquad \mathcal{F}_i^b(t) := \{\, h_j(t) \geq h_j(0) + 1 \text{ for all } j \in \{i, \ldots, \lfloor bi \rfloor\}\}$

for a real $b \in [0, 1)$.

LEMMA 5.4. *Let $\widetilde{C} : \mathbb{R}_+ \to \mathbb{R}$ be a function such that $\widetilde{C}(K) \to \infty$ as $K \to \infty$. Then for each $b \in [0, 1)$, there exists another function $C_b : \mathbb{R}_+ \to \mathbb{R}$ such that $C_b(t) \to \infty$ as $t \to 0$ and the following claim holds.*

*Let $-\infty < \ell < 0 < \mathfrak{r} < \infty$, and let $\underline{\omega}(t)$ be either an $[\ell, \mathfrak{r}]$-monotone process or an $(\ell, \mathfrak{r}, \theta)$-process. Fix a set $U \subseteq \mathbb{Z}$ and a time $t \in (0, \infty)$. Assume that for all $K > 0$ and $i \in U$,*

(5.7) $\qquad \mathbf{P}\left\{\frac{1}{|i|} \sum_{j=i+1}^{0} |\omega_j(t)| \geq K\right\} \leq e^{-\widetilde{C}(K)|i|}.$

*Then for $b \in [0, 1)$, there exists an integer $i_0 = i_0(b) < 0$ such that, for all $i < i_0$ in $U$,*

(5.8) $\qquad \mathbf{P}\{\mathcal{F}_i^b(t)\} \leq e^{-C_b(t)|i|}.$

*The corresponding statement holds under the corresponding conditions for $i > 0$. The conclusion is independent of $\ell, \mathfrak{r}$. Moreover, the statement is also true for the limit of the monotone processes, provided condition (5.7) holds uniformly for the $[\ell, \mathfrak{r}]$-monotone processes.*

PROOF. For this proof always $i \in U$. Consider first the $[\ell, \mathfrak{r}]$-monotone or the $(\ell, \mathfrak{r}, \theta)$-process. For a given $K$, define the event

$$\mathcal{K}_i(t) := \left\{\frac{1}{|i|} \sum_{j=i+1}^{0} |\omega_j(t)| < K\right\}.$$

By assumption (5.7),

$$\mathbf{P}\{\mathcal{F}_i^b(t)\} \leq \mathbf{P}\{\mathcal{F}_i^b(t) \cap \mathcal{K}_i(t)\} + e^{-\widetilde{C}(K)|i|}.$$

On the event $\mathcal{F}_i^b(t)$ each column indexed $i, i+1, \ldots, \lfloor bi \rfloor$ grew by at least one during $[0, t]$. Let $s_1$ denote the time when column $j$ last grew, and $s_2$ the last growth time of that column before $s_1$. If there is no such time, then $s_2 = 0$. By the monotonicity of columns and the rate function $r$, for any $s_2 < s < s_1$,

$$\omega_j(s) = h_{j-1}(s) - h_j(s) \leq h_{j-1}(t) - [h_j(t) - 1] = \omega_j(t) + 1,$$

(5.9) $\quad -\omega_{j+1}(s) = h_{j+1}(s) - h_j(s) \leq h_{j+1}(t) - [h_j(t) - 1] = -\omega_{j+1}(t) + 1$

and

$$r[\omega_j(s)] + r[-\omega_{j+1}(s)] \leq r[\omega_j(t) + 1] + r[-\omega_{j+1}(t) + 1].$$

On the event $\mathcal{K}_i(t)$ we can apply the previous lemma to $\underline{\omega}(t)$. Fix $0 < \alpha < 1 - b$. Then as long as

$$z \geq \frac{3K}{1 - b - \alpha} \quad \text{and} \quad i < i_0(b, \alpha),$$

on the event $\{\mathcal{F}_i^b(t) \cap \mathcal{K}_i(t)\}$ at least $\lfloor \alpha |i| \rfloor$ indices $j$ lie in $\mathcal{D}(z)$ defined by (5.3). For any such index $j$, we have by (5.9) this uniform bound on the rate of column $h_j$:

$$r[\omega_j(s)] + r[-\omega_{j+1}(s)] \leq r[\omega_j(t) + 1] + r[-\omega_{j+1}(t) + 1] < r(z+1) + r(z+1).$$

In other words, in the time interval $[0, t]$, at least $\lfloor \alpha |i| \rfloor$ columns grew during a period when the rate was bounded by $2r(z+1)$. This is true whenever $i < i_0 = i_0(\alpha, b)$. There is no dependence on $\ell$ here, for indices $j < \ell$ do not present a problem. In the $[\ell, \mathfrak{r}]$-monotone process columns outside $\{\ell, \ldots, \mathfrak{r} - 1\}$ do not grow, while for the $(\ell, \mathfrak{r}, \theta)$-process, the leftmost and rightmost constant rates $e^\theta$ and $e^{-\theta}$ are also bounded by $r(z+1)$ if we take $z$ is large enough.

We turn these considerations into a bound:

$\mathbf{P}\{\mathcal{F}_i^b(t) \cap \mathcal{K}_i(t)\}$

$\quad \leq \mathbf{P}\{\text{at least } \lfloor \alpha |i| \rfloor \text{ columns in } i, \ldots, \lfloor bi \rfloor \text{ grew during a time}$

$\quad\quad\quad\quad\quad\quad\quad\quad\quad \text{period when the rate } \leq 2r(z+1)\}$

$\quad \leq \sum_{\substack{\mathcal{D} \subset \{i, \ldots, \lfloor bi \rfloor\} \\ |\mathcal{D}| = \lfloor \alpha |i| \rfloor}} \mathbf{P}\{\text{columns in } \mathcal{D} \text{ grew during a time}$

$\quad\quad\quad\quad\quad\quad\quad\quad\quad \text{period when the rate } \leq 2r(z+1)\}$

$\quad \leq \binom{\lfloor bi \rfloor - i + 1}{\lfloor \alpha |i| \rfloor} \cdot (1 - e^{-2tr(z+1)})^{\lfloor \alpha |i| \rfloor}$

$\quad \leq 2^{\lfloor bi \rfloor + |i| + 1} \cdot (1 - e^{-2tr(z+1)})^{\lfloor \alpha |i| \rfloor}.$

Therefore,

$$\mathbf{P}\{\mathcal{F}_i^b(t)\} \leq 2^{\lfloor bi \rfloor + |i| + 1} \cdot (1 - e^{-2tr(z+1)})^{\lfloor \alpha |i| \rfloor} + e^{-\widetilde{C}(K)|i|}$$

whenever $i < i_0 = i_0(\alpha, b)$ and $z \geq \frac{3K}{1-b-\alpha}$. Move $i$ further to $i < \min\{-2/\alpha, i_0(\alpha, b)\}$ so that $\lfloor \alpha |i| \rfloor / |i| > \alpha - 1/|i| > \alpha/2$. Then for all $K > 0$ and $z \geq \frac{3K}{1-b-\alpha}$,

$$\mathbf{P}\{\mathcal{F}_i^b(t)\} \leq 2 \cdot [2 \cdot (1 - e^{-2tr(z+1)})^{\alpha/2} + e^{-\widetilde{C}(K)}]^{|i|}.$$

The parameter $\alpha$ is fixed for us, so define $i_0(b) = \min\{-2/\alpha, i_0(\alpha, b)\}$. From the bound above we indicate how the function $C_b$ is constructed. Let $A > 0$ be arbitrary:

- Fix $K = K(A)$ large enough that $e^{-\widetilde{C}(K)} < \frac{1}{4}e^{-A}$.
- Fix $z = z(A) \geq \frac{3K}{1-b-\alpha}$.
- Define $\delta = \delta(A) > 0$ by the equality

$$2 \cdot (1 - e^{-2\delta r(z+1)})^{\alpha/2} = \tfrac{1}{4}e^{-A}.$$

This defines a function $C_b(\delta) := A(\delta) \nearrow \infty$ as $\delta \searrow 0$. [Actually, $C_b(\delta)$ is now defined for some finite range $\delta \in (0, \delta_0)$. Put $C_b(\delta) = 0$ for $\delta \geq \delta_0$.]

After these choices, if $t = \delta \in (0, \delta_0)$, the earlier bound turns into (5.8):

$$\mathbf{P}\{\mathcal{F}_i^b(t)\} \leq 2 \cdot (\tfrac{1}{4}e^{-A} + \tfrac{1}{4}e^{-A})^{|i|} \leq e^{-C_b(t)|i|}.$$

For $t \geq \delta_0$, (5.8) is trivially true because $C_b(t) = 0$.

The extension from the $[\ell, \mathfrak{r}]$-monotone process to the case of the limiting process follows by monotonicity of the event $\mathcal{F}_i^b(t)$ in $\ell, \mathfrak{r}$ and the fact that our bounds were independent of these parameters. □

COROLLARY 5.5. *Let $\underline{\zeta}(t)$ be the $(\ell, \mathfrak{r}, \theta_1)$-process, started in distribution $\underline{\pi}$ (4.5), $0 \leq b < 1$ a real number, and*

(5.10) $\mathcal{E}_i^b(t) := \{$*each column $g_j$ with $i \leq j \leq \lfloor bi \rfloor$ of $\underline{\zeta}$ has grown by at least one by time $t$*$\}.$

*Then there is an integer $i_0 = i_0(b)$ and a function $C(t) = C_{b,\theta_1,\theta_2}(t)$, both independent of $\ell, \mathfrak{r}$, such that $C(t) \to \infty$ as $t \to 0$ and*

$$\mathbf{P}\{\mathcal{E}_i^b(t)\} \leq e^{-C(t)|i|} \qquad \text{for all } i < i_0.$$

*The corresponding statement holds under the corresponding conditions for $i > 0$.*

PROOF. Whenever $\ell - 1 \leq i \leq \mathfrak{r}$, Lemma 5.1 and, hence, Lemma 5.4 apply. For $i$ outside of this interval, the probability above is simply zero. $\square$

We now bound exponential moments of column $g_0$ in the $\underline{\zeta}$ process. We let $-\mathcal{L} = \mathcal{R} > 0$ and take $\underline{\zeta}$ as an $(\mathcal{L}, \mathcal{R}, \theta_1)$-process. We first derive a naive bound that depends on the size $\mathcal{L}$ of the process, and then improve it to an estimate independent of $\mathcal{L}$. In an $(\mathcal{L}, \mathcal{R}, \theta_1)$-process $g_{\mathcal{L}-1}$ is the leftmost and $g_\mathcal{R}$ the rightmost column to grow. The boundary increments $\zeta_{\mathcal{L}-1}$ and $\zeta_{\mathcal{R}+1}$ are outside the spatial range that couples with stationary processes (the coupling of Lemma 4.1). However, we can still take advantage of the product-form *initial* distribution.

LEMMA 5.6. *Let $\underline{\zeta}(t)$ be an $(\mathcal{L}, \mathcal{R}, \theta_1)$-process with $\mathcal{L} = -\mathcal{R}$ started from initial distribution $\underline{\pi}$ of* (4.5), *with $g_0(0) = 0$. Then there are finite constants $A$ and $B$ that depend on $K$, $\theta_1$ and $\theta_2$, but not on $\mathcal{L}$ and $t$, such that, for any $K > 0$, $\mathcal{L} < 0$ and $t \geq 0$,*

$$\mathbf{E}(e^{K g_0(t)}) \leq e^{|\mathcal{L}|(At+B)}.$$

PROOF. Utilizing increments as in (2.4) with $g_0(0) = 0$ gives

$$g_0(t) \leq \max_{j=\mathcal{L}-1,\ldots,\mathcal{R}} g_j(t) \leq \max_{j=\mathcal{L}-1,\ldots,\mathcal{R}} g_j(t) - \max_{j=\mathcal{L}-1,\ldots,\mathcal{R}} g_j(0) + \sum_{j=\mathcal{L}}^{\mathcal{R}} |\zeta_j(0)|.$$

Next apply Hölder's inequality:

$$\mathbf{E}(e^{K g_0(t)}) \leq [\mathbf{E}(e^{2K(\max_{j=\mathcal{L}-1,\ldots,\mathcal{R}} g_j(t) - \max_{j=\mathcal{L}-1,\ldots,\mathcal{R}} g_j(0))}) \cdot \mathbf{E}(e^{2K \sum_{j=\mathcal{L}}^{\mathcal{R}} |\zeta_j(0)|})]^{1/2}$$
$$\leq e^{[(4|\mathcal{L}|+2)r(0) + e^{\theta_1} + e^{-\theta_1}][\exp(2K)-1]t} \cdot e^{C'|\mathcal{L}|}.$$

The last inequality comes from the Poisson bound on the growth of the maximum explained above Proposition 3.1, and from the initial product distribution. $\square$

LEMMA 5.7. *Let $\underline{\zeta}(t)$ be an $(\mathcal{L}, \mathcal{R}, \theta_1)$-process with $\mathcal{L} = -\mathcal{R}$, started from initial distribution $\underline{\pi}$ of* (4.5), *with $g_0(0) = 0$. Then for $K > 0$, there exists a finite strictly positive time $T = T(K, \theta_1, \theta_2)$ such that*

$$\mathbf{E}(e^{K g_0(t)})$$

*is bounded uniformly in $\mathcal{L} < 0$ and $0 \leq t < T$.*

PROOF. Fix $b = 0$, and recall the definition (5.10) of $\mathcal{E}_i^0(t)$ for $i \leq 0$. Let $\mathcal{E}_1^0(t)$ be the probability one event. Define also

$$\widetilde{\mathcal{E}}_i^0(t) := \mathcal{E}_i^0(t) \cap \{\text{column } i-1 \text{ of } \underline{\zeta} \text{ has not grown by time } t\}.$$

Since column $\mathcal{L}-1$ is frozen in this process, precisely one of the events $\widetilde{\mathcal{E}}_1^0(t)$, $\widetilde{\mathcal{E}}_0^0(t)$, ..., $\widetilde{\mathcal{E}}_{\mathcal{L}-1}^0(t)$ happens. Write

$$(5.11) \quad \mathbf{E}(e^{Kg_0(t)}) = \sum_{i=\mathcal{L}}^{1} \mathbf{E}(e^{Kg_0(t)} \cdot \mathbf{1}\{\widetilde{\mathcal{E}}_i^0(t)\}) + \mathbf{E}(e^{Kg_0(t)} \cdot \mathbf{1}\{\widetilde{\mathcal{E}}_{\mathcal{L}-1}^0(t)\}).$$

For the first term, under the event $\widetilde{\mathcal{E}}_i^0(t)$, we use the following estimate [with the convention $g_0(0) = 0$ and that empty sums are zero]:

$$(5.12) \quad \begin{aligned} g_0(t) &= g_{i-1}(t) - \sum_{j=i}^{0} \zeta_j(t) = g_{i-1}(0) - \sum_{j=i}^{0} \zeta_j(t) \\ &= \sum_{j=i}^{0} \zeta_j(0) - \sum_{j=i}^{0} \zeta_j(t) \le \sum_{j=i}^{0} \zeta_j(0) - \sum_{j=i}^{0} \eta_j(t). \end{aligned}$$

Hence, with the convention that empty products have value one,

$$\begin{aligned} \sum_{i=\mathcal{L}}^{1} \mathbf{E}(e^{Kg_0(t)} \cdot \mathbf{1}\{\widetilde{\mathcal{E}}_i^0(t)\}) &\le \sum_{i=\mathcal{L}}^{1} \mathbf{E}\left(\prod_{j=i}^{0} e^{K\zeta_j(0)} \prod_{j=i}^{0} e^{-K\eta_j(t)} \cdot \mathbf{1}\{\widetilde{\mathcal{E}}_i^0(t)\}\right) \\ &\le \sum_{i=\mathcal{L}}^{1} \left[\left(\prod_{j=i}^{0} \mathbf{E}e^{3K\zeta_j(0)}\right)\left(\prod_{j=i}^{0} \mathbf{E}e^{-3K\eta_j(t)}\right) \mathbf{P}\{\mathcal{E}_i^0(t)\}\right]^{1/3} \\ &\le 1 + \sum_{i=\mathcal{L}}^{0} e^{C(|i|+1)} \cdot e^{C(|i|+1)} \cdot [\mathbf{P}\{\mathcal{E}_i^0(t)\}]^{1/3}. \end{aligned}$$

The second inequality comes from Hölder and the independence of increments in the configurations $\underline{\eta}(t)$ and $\underline{\zeta}(0)$. The last inequality comes from finite exponential moments of equilibrium increments [see (A.1)] with $C = C(K)$. We apply Corollary 5.5. Take $i_0$ corresponding to $b = 0$ in Corollary 5.5, and split the summation as

$$\sum_{i=\mathcal{L}}^{1} \mathbf{E}(e^{Kg_0(t)} \cdot \mathbf{1}\{\widetilde{\mathcal{E}}_i^0(t)\})$$

$$\le 1 + \sum_{i=\mathcal{L}}^{i_0-1} e^{2C(|i|+1)} \cdot [\mathbf{P}\{\mathcal{E}_i^0(t)\}]^{1/3} + \sum_{i=i_0 \vee \mathcal{L}}^{0} e^{2C(|i|+1)} \cdot [\mathbf{P}\{\mathcal{E}_i^0(t)\}]^{1/3}.$$

The first and the last term on the right are bounded uniformly in $\mathcal{L}$. Apply Corollary 5.5 to the second term to write

$$\sum_{i=\mathcal{L}}^{i_0-1} e^{2C(|i|+1)} \cdot [\mathbf{P}\{\mathcal{E}_i^0(t)\}]^{1/3} \le \sum_{i=\mathcal{L}}^{i_0-1} e^{2C(|i|+1)} \cdot e^{-C(t)|i|/3},$$

which is uniformly bounded in $\mathcal{L}$ if $t$ is taken small enough to guarantee $6C < C(t)$. The function $C(t)$ only depends on $b$, which is 0 now, and $\theta_1$ and $\theta_2$. Hence, the time $T$ only depends on $(K, \theta_1, \theta_2)$.

The second term of (5.11) cannot be handled in a similar fashion since an estimate using (5.12) under the event $\widetilde{\mathcal{E}}^0_{\mathcal{L}-1}(t)$ would involve the variable $\zeta_{\mathcal{L}-1}(t)$ for which we have no convenient bound. Here we use the naive estimate of Lemma 5.6:

$$\mathbf{E}(e^{Kg_0(t)} \cdot \mathbf{1}\{\widetilde{\mathcal{E}}^0_{\mathcal{L}-1}(t)\}) \leq \mathbf{E}(e^{Kg_0(t)} \cdot \mathbf{1}\{\mathcal{E}^0_{\mathcal{L}-1}(t)\})$$
$$\leq [\mathbf{E}(e^{2Kg_0(t)}) \cdot \mathbf{P}\{\mathcal{E}^0_{\mathcal{L}-1}(t)\}]^{1/2}$$
$$\leq e^{|\mathcal{L}|(At+B)} e^{-C(t)|\mathcal{L}-1|/2},$$

with constants $A$, $B$, and a function $C(t) \to \infty$ as $t \to 0$, each also depending on $K$, $\theta_1$, $\theta_2$ by Corollary 5.5 and Lemma 5.6. There is a $T$, depending on $K$, $\theta_1$ and $\theta_2$, such that, for each $0 \leq t \leq T$, $2At + 2B < C(t)$. For such $t$'s, the above expression is bounded from above uniformly in $\mathcal{L}$, which finishes the proof. $\square$

COROLLARY 5.8. *Let $\underline{\zeta}(t)$ be an $(\mathcal{L}, \mathcal{R}, \theta_1)$-process with $\mathcal{L} = -\mathcal{R}$, initially distributed according to $\underline{\pi}$ of (4.5), with $g_0(0) = 0$. Then for any $K > 0$, there exist a time $T > 0$ and constants $\widetilde{D}$ and $\widetilde{E}$, each also depending on $K$, $\theta_1$ and $\theta_2$ but not on $\mathcal{L}$, such that*

$$\mathbf{E}(e^{Kg_i(t)}) \leq e^{\widetilde{D} + \widetilde{E}|i|}$$

*for any $\mathcal{L} - 1 \leq i \leq \mathcal{R}$ and $0 \leq t < T$.*

PROOF. For $i < 0$, (2.4) reads as

$$g_i(t) = \sum_{j=i+1}^{0} \zeta_j(t) + g_0(t) \leq \sum_{j=i+1}^{0} \xi_j(t) + g_0(t)$$

by the coupling of $\underline{\zeta}$ and $\underline{\xi}$. Hence, by Hölder, we have

$$\mathbf{E}(e^{Kg_i(t)}) \leq [\mathbf{E}(e^{2K\sum_{j=i+1}^{0} \xi_j(t)}) \mathbf{E}(e^{2Kg_0(t)})]^{1/2} \leq e^{\widetilde{E}|i|} \cdot e^{\widetilde{D}}$$

up to some time $T$, with some constants $\widetilde{E}$ and $\widetilde{D}$ by the product equilibrium of $\underline{\xi}(t)$, by (A.1), and Lemma 5.7. A similar argument comparing $-\zeta_j$'s with $-\eta_j$'s proves the statement for positive $i$'s. $\square$

5.3. *Further regularity properties of the $\underline{\omega}(t)$ process.* Now we refine Lemma 5.2 to get exponential moment bounds on sums of $\omega_j(t)$-variables. Consequences include moment bounds on $\underline{\omega}$-rates and bounds on growth of adjacent columns.

Note that while Lemma 5.2 was valid for all time, the exponential moments we proved for $\underline{\zeta}$ and use below are valid only up to some time $T = T(\theta_1, \theta_2) > 0$. We still use the conditional coupling of Section 4.2. Recall the definition (2.8) of $K_{\underline{\omega}}$. The main point is again to get estimates independent of $\ell, \mathfrak{r}$ so that we can go to the limit and get estimates for the infinite-volume processes.

LEMMA 5.9. *Let $\underline{\omega}(t)$ be either an $[\ell, \mathfrak{r}]$-monotone process or the limit of these processes started from an initial state $\underline{\omega} \in \widetilde{\Omega}$. Let $K > 0$. Then there is a finite constant $C = C(K, \underline{\omega})$ and a positive time $T = T(K, K_{\underline{\omega}})$ such that, for any $t < T$,*

$$\mathbf{E}^{\underline{\omega}}(e^{K \sum_{j=i+1}^{0} |\omega_j(t)|}) \leq e^{C|i|}$$

*holds whenever $\ell \leq i < 0$, and*

$$\mathbf{E}^{\underline{\omega}}(e^{K \sum_{j=1}^{i} |\omega_j(t)|}) \leq e^{Ci}$$

*holds whenever $1 \leq i < \mathfrak{r}$.*

PROOF. Consider first the $[\ell, \mathfrak{r}]$-monotone process $\underline{\omega}^{[\ell, \mathfrak{r}]}(t)$. Fix $\mathcal{L} = -\mathcal{R}$ such that $\mathcal{L} \leq \ell < \mathfrak{r} \leq \mathcal{R}$. In the conditional coupling argument, we can couple $\underline{\omega}^{[\ell, \mathfrak{r}]}(t)$ with the $(\mathcal{L}, \mathcal{R}, \theta_1)$-process $\underline{\zeta}$. The arguments of the proof of Lemma 5.2 remain valid with this setting. Take (5.2) from that proof with $t_0 = t = t_1$. We show that each term of the right-hand side of (5.2) has the stated exponential moment, and then Hölder's inequality leads to the proof.

The first term with the sum of terms $K|\omega_j(0)|$ is simply smaller than $\widetilde{K}|i|$ for some $\widetilde{K}$ since $\underline{\omega} \in \widetilde{\Omega}$.

For the second and third terms, we write, both for time 0 and time $t$,

$$\mathbf{E}(e^{K \sum_{j=i+1}^{0} |\zeta_j|}) \leq \mathbf{E}(e^{K \sum_{j=i+1}^{0} (|\eta_j| + |\xi_j|)})$$
$$\leq [\mathbf{E}(e^{2K \sum_{j=i+1}^{0} |\eta_j|}) \mathbf{E}(e^{2K \sum_{j=i+1}^{0} |\xi_j|})]^{1/2} \leq e^{C|i|}$$

by the coupling $\eta_j \leq \zeta_j \leq \xi_j$ and the product equilibrium of $\underline{\eta}$ and $\underline{\xi}$.

The fourth term is now zero.

For the term $h_i^{[\ell, \mathfrak{r}]}(t)$ we again use conditional coupling:

$$\mathbf{E}^{\underline{\omega}}(e^{K h_i^{[\ell, \mathfrak{r}]}(t)}) = \mathbf{E}^{(\underline{\pi}|\mathcal{A}^{\underline{\omega}})}(e^{K h_i^{[\ell, \mathfrak{r}]}(t)})$$

(5.13)

$$\leq \mathbf{E}^{(\underline{\pi}|\mathcal{A}^{\underline{\omega}})}(e^{K g_i^{(\mathcal{L}, \mathcal{R}, \theta_1)}(t)}) \leq \frac{\mathbf{E}(e^{K g_i^{(\mathcal{L}, \mathcal{R}, \theta_1)}(t)})}{\underline{\pi}\{\underline{\zeta} \in \mathcal{A}^{\underline{\omega}}\}},$$

and a similar argument applies to the last term $h_0^{[\ell,\mathfrak{r}]}(t)$.

All that is left is to bound the exponential moments of $g_i^{(\mathcal{L},\mathcal{R},\theta_1)}(t)$ and $g_0^{(\mathcal{L},\mathcal{R},\theta_1)}(t)$. This is done by Corollary 5.8 up to a finite time $T = T(K, \theta_1, \theta_2)$. The parameters $\theta_1$ and $\theta_2$ are chosen according to (4.4), hence, $T$ acquires a dependence on $K_{\underline{\omega}}$.

Hölder's inequality now together with (5.2) proves the statement for $i < 0$ and finite $\ell$, $\mathfrak{r}$. A similar argument works for $i > 0$. Notice that none of our bounds depends on $\ell$ and $\mathfrak{r}$ as long as $\ell \leq i < \mathfrak{r}$ holds. Fatou's lemma therefore extends the result for the limit $\ell \to -\infty$, $\mathfrak{r} \to \infty$. $\square$

From these moment bounds we get hypothesis (5.7) for controlling the growth of a large number of adjacent columns. Recall definition (5.6).

COROLLARY 5.10. *Let $\underline{\omega}(t)$ be either an $[\ell, \mathfrak{r}]$-monotone process or the limit of these processes started from an initial state $\underline{\omega} \in \widetilde{\Omega}$. Let $b \in [0, 1)$. Then there exist a time $T^{\underline{\omega}} > 0$, an integer $i_0 = i_0(b) < 0$ and a function $C_b^{\underline{\omega}} : \mathbb{R}^+ \to \mathbb{R}$ with these properties: $C_b^{\underline{\omega}}(t) \to \infty$ as $t \to 0$, and for all $i < i_0$,*

$$(5.14) \qquad \mathbf{P}^{\underline{\omega}}\{\mathcal{F}_i^b(t)\} \leq e^{-C_b^{\underline{\omega}}(t)|i|}$$

*whenever $t < T^{\underline{\omega}}$. The corresponding statement holds under the corresponding conditions for $i > 0$. Note again the independence of the estimates of $\ell, \mathfrak{r}$.*

PROOF. We need to verify that (5.7) holds. For $\ell \leq i < 0$ and $t < T$, we write

$$\mathbf{P}^{\underline{\omega}}\left\{\frac{1}{|i|}\sum_{j=i+1}^{0}|\omega_j(t)| \geq K\right\} \leq e^{-K|i|}\mathbf{E}^{\underline{\omega}}(e^{\sum_{j=i+1}^{0}|\omega_j(t)|}) \leq e^{(C-K)|i|}$$

by Markov's inequality and Lemma 5.9. For $i < \ell$, $\mathbf{P}^{\underline{\omega}}\{\mathcal{F}_i^b(t)\}$ is simply zero. A similar argument holds for positive $i$'s. Notice that $C$ on the right-hand side above depends on $\underline{\pi}\{\underline{\zeta} \in \mathcal{A}^{\underline{\omega}}\}$ of the conditional coupling used in Lemma 5.9. This dependence is then transferred to $\widetilde{C}(K)$ in (5.7), and from there to $C_b(t)$. $\square$

We derive one last estimate for the next section.

COROLLARY 5.11. *Under the hypotheses of Lemma 5.9, for any $1 \leq p < \infty$, there exist $T = T(p, \underline{\omega}) > 0$, $A = A(p, \underline{\omega})$ and $B = B(p, \underline{\omega})$ constants such that*

$$\mathbf{E}^{\underline{\omega}}\left[\sup_{0 \leq t \leq T} r(|\omega_i(t)|)^p\right] \leq e^{A+B|i|}.$$

PROOF. By the exponential bound (2.3) on the rates and by monotonicity of columns, for any $0 \le t \le T$,

$$r(|\omega_i(t)|)^p \le e^{\beta p |\omega_i(t)|} = e^{\beta p |h_{i-1}(t) - h_i(t)|}$$
$$\le e^{\beta p (|h_{i-1}(0)| - h_{i-1}(0) + |h_i(0)| - h_i(0))} \cdot e^{\beta p h_{i-1}(T)} \cdot e^{\beta p h_i(T)}.$$

Take expectations, apply the Schwarz inequality, apply the conditional coupling as in (5.13), then Corollary 5.8, and finally let $-\ell, \mathfrak{r} \to \infty$. □

**6. Analytic properties.** Fix a bounded cylinder function $\varphi$ on $\widetilde{\Omega}$. For each finite-volume process $\underline{\omega}^{[\ell, \mathfrak{r}]}(\cdot)$, we have the integrated forward and backward equations

$$(6.1) \quad \begin{aligned} S^{[\ell, \mathfrak{r}]}(t) \varphi(\underline{\omega}) - \varphi(\underline{\omega}) &= \int_0^t S^{[\ell, \mathfrak{r}]}(s) L^{[\ell, \mathfrak{r}]} \varphi(\underline{\omega}) \, ds \\ &= \int_0^t L^{[\ell, \mathfrak{r}]}(S^{[\ell, \mathfrak{r}]}(s) \varphi)(\underline{\omega}) \, ds, \end{aligned}$$

where $S^{[\ell, \mathfrak{r}]}(t)$ is the semigroup of the process. The generator $L^{[\ell, \mathfrak{r}]}$ may involve exponentials of the unbounded variables $(\omega_i^{[\ell, \mathfrak{r}]})$, so an argument is needed to justify the existence of the right-hand side of the first line of (6.1). But it can be checked that the process

$$t \mapsto \sum_{i=\ell}^{\mathfrak{r}} |\omega_i^{[\ell, \mathfrak{r}]}(t)| - \sum_{i=\ell}^{\mathfrak{r}} |\omega_i^{[\ell, \mathfrak{r}]}(0)|$$

is bounded by twice a Poisson process of rate $2r(0)(\ell - \mathfrak{r})$.

In this section we seek to let $-\ell, \mathfrak{r} \to \infty$ in (6.1) to capture corresponding equations for the limiting process $\underline{\omega}(\cdot)$. The left-hand side of (6.1) converges to $S(t)\varphi(\underline{\omega})$ by the definition of the process as a limit. The right-hand sides need some work and restrictive assumptions. Throughout the section the initial state satisfies $\underline{\omega} \in \widetilde{\Omega}$.

6.1. *Integrated forward equation.* Since $\varphi$ is a bounded cylinder function, there is an integer $a < \infty$ such that

$$(6.2) \quad |L\varphi(\underline{\omega}(s))| \le C \sum_{i=-a}^{a} r_i(\underline{\omega}(s)).$$

Together with the moment bound of Corollary 5.11, this suggests that we can verify the equations up to a time $T^{\underline{\omega}}$ that depends on the initial state $\underline{\omega}$. Time integration enables us to do better in the forward equation.

Note that $L^{[\ell,\mathfrak{r}]}\varphi = L\varphi$ as soon as $-\ell, \mathfrak{r}$ are large enough. Write, with (3.3), the right-hand side of the forward equation in (6.1) as

$$\int_0^t S^{[\ell,\mathfrak{r}]}(s) L\varphi(\underline{\omega}) \, ds$$
$$= \sum_{i=-a}^{a} \int_0^t \mathbf{E}^{\underline{\omega}} r_i(\underline{\omega}^{[\ell,\mathfrak{r}]}(s)) [\varphi(\underline{\omega}^{[\ell,\mathfrak{r}](i,i+1)}(s)) - \varphi(\underline{\omega}^{[\ell,\mathfrak{r}]}(s))] \, ds.$$

Fix $i$, and think of the integrand

$$f_{\ell,\mathfrak{r}} = r_i(\underline{\omega}^{[\ell,\mathfrak{r}]}(s))[\varphi(\underline{\omega}^{[\ell,\mathfrak{r}](i,i+1)}(s)) - \varphi(\underline{\omega}^{[\ell,\mathfrak{r}]}(s))]$$

as a function of $(s, \underline{N})$, integrated with respect to $ds \otimes \mathbb{P}$, where $\mathbb{P}$ is the distribution of the family $\underline{N}$ of Poisson processes. [Recall the construction of the processes $\underline{\omega}^{[\ell,\mathfrak{r}]}(\cdot)$ and $\underline{\omega}(\cdot)$ as measurable functions of $\underline{N}$ when the initial state is given.] Almost everywhere convergence as $-\ell, \mathfrak{r} \to \infty$ to the desired limit

$$f = r_i(\underline{\omega}(s))[\varphi(\underline{\omega}^{(i,i+1)}(s)) - \varphi(\underline{\omega}(s))]$$

follows from the construction. We will apply the generalized dominated convergence theorem (e.g., Proposition 18 from Chapter 11 of [13]). The domination is done by

$$|f_{\ell,\mathfrak{r}}| \leq g_{\ell,\mathfrak{r}} := Cr_i(\underline{\omega}^{[\ell,\mathfrak{r}]}(s)) \xrightarrow[-\ell,\mathfrak{r}\to\infty]{} Cr_i(\underline{\omega}(s)) =: g.$$

The key help from the integration of the time variable is that the forward equation of the finite-volume height process $\underline{h}^{[\ell,\mathfrak{r}]}(\cdot)$ gives the bound

$$(6.3) \quad \int_0^t \mathbf{E}^{\underline{\omega}} r_i(\underline{\omega}^{[\ell,\mathfrak{r}]}(s)) \, ds = \mathbf{E}^{\underline{\omega}}[h_i^{[\ell,\mathfrak{r}]}(t) - h_i(0)] \leq A^{\underline{\omega}} t + B^{\underline{\omega}}|i| + |h_i(0)|,$$

where the last inequality came from Theorem 2.4. It only remains to check that

$$(6.4) \qquad \int_0^t \mathbf{E}^{\underline{\omega}} r_i(\underline{\omega}^{[\ell,\mathfrak{r}]}(s)) \, ds \longrightarrow \int_0^t \mathbf{E}^{\underline{\omega}} r_i(\underline{\omega}(s)) \, ds.$$

We take advantage of the freedom in the order of convergence. Lemmas 3.3 and 3.4 gave the monotonicity

$$\omega_i^{[\ell,\mathfrak{r}]}(s) \leq \omega_i^{[\ell-1,\mathfrak{r}]}(s) \quad \text{and} \quad \omega_i^{[\ell,\mathfrak{r}]}(s) \geq \omega_i^{[\ell,\mathfrak{r}+1]}(s)$$

(6.5)
$$\text{for any } \ell \leq i \leq \mathfrak{r}.$$

By (3.3), the integrands are sums of two terms of the form $r(\cdot)$ for all large $\mathfrak{r}$ and large negative $\ell$. Monotonicity of the rate function $r$ gives similar inequalities for these terms separately, we show how the limit (6.4) is done for the first term:

$$\int_0^t \mathbf{E}^{\underline{\omega}} r(\underline{\omega}_i^{[\ell,\mathfrak{r}]}(s))\,ds \xrightarrow[\mathfrak{r}\to\infty]{} \int_0^t \mathbf{E}^{\underline{\omega}} r(\underline{\omega}_i^{[\ell,\infty]}(s))\,ds \xrightarrow[\ell\to-\infty]{} \int_0^t \mathbf{E}^{\underline{\omega}} r(\underline{\omega}_i(s))\,ds.$$

The first limit uses the dominated convergence theorem with the bound (6.3) for each fixed $\ell$. The limit $\underline{\omega}^{[\ell,\infty]}(s)$ is by definition the configuration of increments of the height configuration $\underline{h}^{[\ell,\infty]}(s)$ which itself is the well-defined and finite monotone limit of $\underline{h}^{[\ell,\mathfrak{r}]}(s)$ as $\mathfrak{r}\to\infty$. The variable $\omega_i^{[\ell,\mathfrak{r}]}(s)$ does not change once $\mathfrak{r}$ is large enough, and so the other monotonicity from (6.5) becomes $\omega_i^{[\ell,\infty]}(s) \leq \omega_i^{[\ell-1,\infty]}(s)$. Now the second limit above comes from monotone convergence. The limit of $r(-\omega_{i+1}^{[\ell,\mathfrak{r}]}(s))$ is treated similarly.

We have proved the convergence

$$\int_0^t S^{[\ell,\mathfrak{r}]}(s) L\varphi(\underline{\omega})\,ds$$

$$= \sum_{i=-a}^{a} \int_0^t \mathbf{E}^{\underline{\omega}} r_i(\underline{\omega}^{[\ell,\mathfrak{r}]}(s))[\varphi(\underline{\omega}^{[\ell,\mathfrak{r}](i,i+1)}(s)) - \varphi(\underline{\omega}^{[\ell,\mathfrak{r}]}(s))]\,ds$$

$$\longrightarrow \sum_{i=-a}^{a} \int_0^t \mathbf{E}^{\underline{\omega}} r_i(\underline{\omega}(s))[\varphi(\underline{\omega}^{(i,i+1)}(s)) - \varphi(\underline{\omega}(s))]\,ds$$

$$= \int_0^t S(s) L\varphi(\underline{\omega})\,ds.$$

The integrability of $L\varphi(\underline{\omega}(s))$ over $[0,T] \times D_{\widetilde{\Omega}}$ is contained in the estimates above: apply Fatou's lemma to the left-hand side of (6.3) to get an inequality that shows the integrability of the upper bound in (6.2). This completes the proof of Theorem 2.6.

6.2. *Integrated backward equation.* To see that $S(s)\varphi(\underline{\omega})$ is a function to which the infinite generator $L$ can be applied, we need some bounds.

LEMMA 6.1. *There is a function $A(s,\underline{\omega}) \nearrow \infty$ as $s \searrow 0$ that depends on the function $\varphi$, but not on $[\ell,\mathfrak{r}]$, such that*

$$(6.6) \qquad |S^{[\ell,\mathfrak{r}]}(s)\varphi(\underline{\omega}^{(j,j+1)}) - S^{[\ell,\mathfrak{r}]}(s)\varphi(\underline{\omega})| \leq Ce^{-A(s,\underline{\omega})|j|}$$

*for all $|j| \geq j_0$. $C$ and $j_0$ are fixed constants independent of $\underline{\omega}$, $\ell$ and $\mathfrak{r}$. This bound also holds in the limit $\ell \to -\infty$, $\mathfrak{r} \to \infty$ for the limiting process.*

PROOF. We consider $j < 0$ and omit the symmetric argument for $j > 0$. We couple two $[\ell,\mathfrak{r}]$ processes $\underline{\omega}'(\cdot)$ and $\underline{\omega}(\cdot)$ started at $\underline{\omega}'(0) = \underline{\omega}^{(j,j+1)}$ and $\underline{\omega}(0) = \underline{\omega}$. The difference between the processes can be represented by a second class antiparticle $Q_\downarrow$ and a second class particle $Q_\uparrow$ added to the process $\underline{\omega}(\cdot)$ at time zero. The $(Q_\downarrow, Q_\uparrow)$ pair annihilates each other if they meet. After this they are no longer in the system and we set $Q_\downarrow = -\infty$ and $Q_\uparrow = \infty$. No other second class particles are created in the process.

As before, $[-a, a]$ represents a bounded interval of sites that determine the value $\varphi(\underline{\omega})$. Suppose $\ell - 1 \leq j < -2a < -6$. (If $j \leq \ell - 2$, the extra brick over $[j, j+1]$ does not affect an $[\ell, \mathfrak{r}]$-process.) Since $\varphi$ sees the difference between $\underline{\omega}'(s)$ and $\underline{\omega}(s)$ only if it sees a second class particle,

$$
\begin{aligned}
|S^{[\ell,\mathfrak{r}]}&(s)\varphi(\underline{\omega}^{(j,j+1)}) - S^{[\ell,\mathfrak{r}]}(s)\varphi(\underline{\omega})| \\
&\leq \mathbf{E}[|\varphi(\underline{\omega}'^{[\ell,\mathfrak{r}]}(s)) - \varphi(\underline{\omega}^{[\ell,\mathfrak{r}]}(s))||\underline{\omega}'^{[\ell,\mathfrak{r}]}(0) = \underline{\omega}^{(j,j+1)}, \underline{\omega}^{[\ell,\mathfrak{r}]}(0) = \underline{\omega}] \\
&\leq C \mathbf{P}^{\underline{\omega}} \{ -a \leq Q_\downarrow(s) \leq a \text{ or } -a \leq Q_\uparrow(s) \leq a | Q_\downarrow(0) = j, Q_\uparrow(0) = j+1 \} \\
&\leq C \mathbf{P}^{\underline{\omega}} \{ -a \leq Q_\uparrow(s) < \infty | Q_\downarrow(0) = j, Q_\uparrow(0) = j+1 \}.
\end{aligned}
$$

The last event above implies that all sites in $\{j+1, \ldots, -a-1\}$ have experienced a jump by time $s$. The coupling keeps the columns of the process $\underline{\omega}'(\cdot)$ at least as high as those of the process $\underline{\omega}(\cdot)$. Consequently, process $\underline{\omega}'(\cdot)$ has experienced event $\mathcal{F}_{j+1}^{3/4}(s)$ of (5.6), and the last line in the calculation above is bounded by $C \mathbf{P}^{\underline{\omega}^{(j,j+1)}} \{ F_{j+1}^{3/4}(s) \}$.

Next we argue that

(6.7)  $\mathbf{P}^{\underline{\omega}^{(j,j+1)}} \{ F_{j+1}^{3/4}(s) \} \leq \exp(-A(s, \underline{\omega})|j|)$

for a function $A(s, \underline{\omega})$ with the properties claimed in the statement of the lemma. Since the system contains at most one second class particle and antiparticle, $|\omega_k(t) - \omega_k'(t)| \leq 1$ for each site $k$ and all time $0 \leq t < \infty$. Therefore, the bounds of Lemma 5.9 hold for $\underline{\omega}'(t)$ if $t \in [0, T(K, K_{\underline{\omega}})]$ and $C(K, \underline{\omega})$ is replaced by $C(K, \underline{\omega}) + K$. The function $C(K, \underline{\omega}) + K$ in turn determines the function that appears in the exponent on the right-hand side of the estimate (5.14). Hence, the statement

$$\mathbf{P}^{\underline{\omega}^{(j,j+1)}} \{ \mathcal{F}_i^b(t) \} \leq e^{-C_b^{\underline{\omega}}(t)|i|}$$

from Corollary 5.10, with an adjusted $C_b^{\underline{\omega}}(t)$, is valid for all $j$ as long as $t < T^{\underline{\omega}}$ and $i < i_0(b)$. Take the case $b = 3/4$, $i = j + 1$, restrict further to $j < i_0(3/4) - 1$, and adjust the exponent slightly to get (6.7). The lemma is proved for negative $j$.  $\square$

Assumption $\underline{\omega} \in \widetilde{\Omega}$ implies the existence of a constant $B(\underline{\omega})$ such that

$$r_i(\underline{\omega}) \leq e^{B(\underline{\omega})|i|} \qquad \text{for all } i.$$

Fix $T^{\underline{\omega}} > 0$ so that $A(t, \underline{\omega}) > B(\underline{\omega}) + 1$ for $0 \leq t \leq T^{\underline{\omega}}$. Then estimate (6.6) guarantees that $L(S(t)\varphi)(\underline{\omega})$ is well defined for $0 \leq t \leq T^{\underline{\omega}}$. We show that the right-hand side of the integrated backward equation for $S^{[\ell,\mathfrak{r}]}(t)$ converges to that of $S(t)$. Let $t \in (0, T^{\underline{\omega}}]$:

$$\text{(6.8)} \quad \int_0^t LS(s)\varphi(\underline{\omega})\,ds - \int_0^t L^{[\ell,\mathfrak{r}]}S^{[\ell,\mathfrak{r}]}(s)\varphi(\underline{\omega})\,ds$$

$$\text{(6.9)} \quad = \int_0^t LS(s)\varphi(\underline{\omega})\,ds - \int_0^t L^{[\ell,\mathfrak{r}]}S(s)\varphi(\underline{\omega})\,ds$$

$$\text{(6.10)} \quad + \int_0^t L^{[\ell,\mathfrak{r}]}S(s)\varphi(\underline{\omega})\,ds - \int_0^t L^{[\ell,\mathfrak{r}]}S^{[\ell,\mathfrak{r}]}(s)\varphi(\underline{\omega})\,ds.$$

Line (6.9) is bounded in absolute value by

$$t \sum_{j \notin [\ell,\mathfrak{r}-1]} r_j(\underline{\omega}) C e^{(-B(\underline{\omega})-1)|j|} \longrightarrow 0 \quad \text{as } -\ell, \mathfrak{r} \to \infty.$$

Line (6.10) equals

$$\sum_{j\in\mathbb{Z}} r_j(\underline{\omega})\mathbf{1}_{\{\ell \le j < \mathfrak{r}\}} \int_0^t (S(s)\varphi(\underline{\omega}^{(j,j+1)}) - S^{[\ell,\mathfrak{r}]}(s)\varphi(\underline{\omega}^{(j,j+1)})$$
$$+ S^{[\ell,\mathfrak{r}]}(s)\varphi(\underline{\omega}) - S(s)\varphi(\underline{\omega}))\,ds.$$

Each term in the sum vanishes as $-\ell, \mathfrak{r} \to \infty$, and (6.6) enables us to apply dominated convergence to the entire sum.

We have shown that (6.8) vanishes in the limit and thereby proved the first statement of Theorem 2.7.

### 6.3. Differentiating the semigroup.

LEMMA 6.2. *There exists a time $T^{\underline{\omega}} > 0$ depending on $\underline{\omega}$ such that $t \mapsto S(t)L\varphi(\underline{\omega})$ and $t \mapsto LS(t)\varphi(\underline{\omega})$ are finite, continuous functions on $[0, T^{\underline{\omega}}]$.*

PROOF. For each finite-volume process and index $j$, $\omega_j^{[\ell,\mathfrak{r}]}(s) \to \omega_j^{[\ell,\mathfrak{r}]}(t)$ a.s. as $s \to t$ because these are well defined jump processes and the probability of a jump at time $t$ is zero. By Proposition 4.5, this convergence $\omega_j(s) \to \omega_j(t)$ as $s \to t$ is valid for the limiting process. Now apply these limits to

$$S(t)L\varphi(\underline{\omega}) = \sum_{i=-a}^{a} \mathbf{E}^{\underline{\omega}} r_i(\underline{\omega}(t))[\varphi(\underline{\omega}^{(i,i+1)}(t)) - \varphi(\underline{\omega}(t))]$$

and recall the moment bounds in Corollary 5.11.

For $LS(t)\varphi(\underline{\omega})$, each term of the sum

$$LS(s)\varphi(\underline{\omega}) = \sum_{i\in\mathbb{Z}} r_j(\underline{\omega})[S(s)\varphi(\underline{\omega}^{(i,i+1)}) - S(s)\varphi(\underline{\omega})]$$

converges to the right limit as $s \to t$. Use estimate (6.6) to apply dominated convergence to the sum as in the previous proof. □

Theorem 2.6 and the part of Theorem 2.7 already proved imply that for $t \in [0, T^{\underline{\omega}}]$,

$$S(t)\varphi(\underline{\omega}) - \varphi(\underline{\omega}) = \int_0^t S(s)L\varphi(\underline{\omega})\,ds = \int_0^t LS(s)\varphi(\underline{\omega})\,ds,$$

and by the lemma above, the integrands are continuous in $s$. This proves the remaining part of Theorem 2.7.

REMARK 6.3. The results of this section can be repeated for a pair of processes, coupled via the basic coupling as described by Tables 3 and 4, provided that both processes are started from a state in $\widetilde{\Omega}$ of (2.9), or from a good distribution of Section 2.3.

**7. Invariance and ergodicity.** We show the invariance of $\underline{\mu}^{(\theta)}$ first. An argument is needed because while Proposition 3.1 showed the invariance of the product measure for an $(\ell, \mathfrak{r}, \theta)$-process (without the boundary increments), the process $\underline{\omega}(\cdot)$ we investigate is a limit of a different finite-volume process, namely, the $[\ell, \mathfrak{r}]$-monotone process. To simplify notation, throughout this section we write $\mathbf{E}$ and $\mathbf{P}$ for $\mathbf{E}^{\underline{\mu}^{(\theta)}}$ and $\mathbf{P}^{\underline{\mu}^{(\theta)}}$.

LEMMA 7.1. *Let $\varphi$ be a bounded cylinder function on $\widetilde{\Omega}$, in other words, $\varphi$ depends on only finitely many $\omega_i$-values. Then*

$$\mathbf{E}[\varphi(\underline{\omega}(t))] = \mathbf{E}[\varphi(\underline{\omega}(0))] \qquad \text{for } 0 \le t < \infty.$$

PROOF. Consider the limiting process $\underline{\omega}(t)$, the $[\ell, \mathfrak{r}]$-monotone process $\underline{\omega}^{[\ell,\mathfrak{r}]}(t)$ and the $(\ell, \mathfrak{r}, \theta)$-process $\underline{\zeta}^{(\ell,\mathfrak{r},\theta)}(t)$, all started from the same initial state $\underline{\omega}(0) = \underline{\zeta}(0) \sim \underline{\mu}^{(\theta)}$. Write

$$(7.1) \quad |\mathbf{E}[\varphi(\underline{\omega}(t))] - \mathbf{E}[\varphi(\underline{\omega}(0))]| \le |\mathbf{E}[\varphi(\underline{\omega}(t))] - \mathbf{E}[\varphi(\underline{\omega}^{[\ell,\mathfrak{r}]}(t))]|$$

$$(7.2) \qquad\qquad\qquad + |\mathbf{E}[\varphi(\underline{\omega}^{[\ell,\mathfrak{r}]}(t))] - \mathbf{E}[\varphi(\underline{\zeta}^{(\ell,\mathfrak{r},\theta)}(t))]|$$

$$(7.3) \qquad\qquad\qquad + |\mathbf{E}[\varphi(\underline{\zeta}^{(\ell,\mathfrak{r},\theta)}(t))] - \mathbf{E}[\varphi(\underline{\zeta}(0))]|.$$

We show that, for some $T = T(\theta) > 0$, the terms on the right-hand side above tend to 0 for each $t \in [0, T]$ as $-\ell, \mathfrak{r} \to \infty$. Invariance in a fixed time interval $[0, T]$ suffices for then the Markov property extends it for all time. □

For the term on the right-hand side of line (7.1), this follows for any $t$ from the construction. Specifically, by Proposition 4.5, the probability that $\underline{\omega}(t)$ and $\underline{\omega}^{[\ell,\mathfrak{r}]}(t)$ differ on the support of $\varphi$ vanishes as $-\ell, \mathfrak{r} \to \infty$. The last line (7.3) vanishes for all $t$ as soon as $-\ell$ and $\mathfrak{r}$ are large enough. The reason is that, except for the boundary increments $\zeta_{\ell-1}^{(\ell,\mathfrak{r},\theta)}(t)$ and $\zeta_{\mathfrak{r}+1}^{(\ell,\mathfrak{r},\theta)}(t)$, the configurations $\underline{\zeta}^{(\ell,\mathfrak{r},\theta)}(t)$ and $\underline{\zeta}(0)$ are equal in distribution.

The term on line (7.2) needs bounds on second class particles. Apply the basic coupling of Section 3.4 between $\underline{\omega}^{[\ell,\mathfrak{r}]}(t)$ and $\underline{\zeta}^{(\ell,\mathfrak{r},\theta)}(t)$. Define the event

$$\mathcal{B}^{[\ell,\mathfrak{r}]}(t) = \{\zeta_j^{(\ell,\mathfrak{r},\theta)}(t) \neq \omega_j^{[\ell,\mathfrak{r}]}(t) \text{ for some } \lfloor \ell/2 \rfloor < j < \lceil \mathfrak{r}/2 \rceil\}.$$

Once the interval $[\ell, \mathfrak{r}]$ is large enough, $\varphi$ will depend only on $\omega$'s over sites in the range $\lfloor \ell/2 \rfloor + 1, \ldots, \lceil \mathfrak{r}/2 \rceil - 1$, and we can write

$$|\mathbf{E}[\varphi(\underline{\omega}^{[\ell,\mathfrak{r}]}(t))] - \mathbf{E}[\varphi(\underline{\zeta}^{(\ell,\mathfrak{r},\theta)}(t))]| \leq C\mathbf{P}\{\mathcal{B}^{[\ell,\mathfrak{r}]}(t)\}.$$

The event $\mathcal{B}^{[\ell,\mathfrak{r}]}(t)$ implies that there is at least one second class particle in $I = \{\lfloor \ell/2 \rfloor + 1, \ldots, \lceil r/2 \rceil - 1\}$. It is a property of the coupling that no second class particles or antiparticles are created in the interior $[\ell+1, \mathfrak{r}-1]$. Second class particles can come into $I$ only from site $\ell$, and second class antiparticles can come into $I$ only from site $\mathfrak{r}$. If a second class particle moves from $\ell$ to $\lfloor \ell/2 \rfloor + 1$ by time $t$, then each intervening site has experienced a jump during $[0,t]$. The coupling also keeps the columns of $\underline{\zeta}^{(\ell,\mathfrak{r},\theta)}$ above those of $\underline{\omega}^{[\ell,\mathfrak{r}]}$ (Lemma 3.2), hence, the columns of $\underline{\zeta}^{(\ell,\mathfrak{r},\theta)}$ must have all increased in the range $\ell, \ldots, \lfloor \ell/2 \rfloor$. In other words, event $\mathcal{E}_\ell^{1/2}(t)$ of (5.10) happened. By the same token, if a second class antiparticle came from $\mathfrak{r}$ to $\lceil \mathfrak{r}/2 \rceil - 1$, event $\mathcal{E}_\mathfrak{r}^{1/2}(t)$ happened.

By Corollary 5.5, $\mathcal{E}_\ell^{1/2}(t)$ and $\mathcal{E}_\mathfrak{r}^{1/2}(t)$ have exponentially small probability in $\ell$ and $\mathfrak{r}$, respectively, until some time $T$ depending on $\theta$. Consequently,

$$\lim_{-\ell,\mathfrak{r}\to\infty} |\mathbf{E}[\varphi(\underline{\omega}^{[\ell,\mathfrak{r}]}(t))] - \mathbf{E}[\varphi(\underline{\zeta}^{(\ell,\mathfrak{r},\theta)}(t))]| = 0$$

for all $t < T$.

We have shown that the terms on lines (7.1)–(7.3) vanish as $-\ell, \mathfrak{r} \to \infty$, and thereby proved the lemma.

The lemma above implies the stationarity, and we turn to ergodicity. First an auxiliary lemma.

LEMMA 7.2.   *Fix* $\underline{\omega} \in \widetilde{\Omega}$ *and* $i \in \mathbb{Z}$. *Define* $\underline{\zeta} := \underline{\omega}^{(i,i+1)}$, *that is,*

(7.4)   $$\zeta_j = \begin{cases} \omega_j, & \text{for } j \neq i, i+1, \\ \omega_j - 1, & \text{for } j = i, \\ \omega_j + 1, & \text{for } j = i+1. \end{cases}$$

*Let the two processes* $\underline{\omega}(t)$ *and* $\underline{\zeta}(t)$ *start from initial states* $\underline{\omega}$ *and* $\underline{\zeta}$, *respectively, and evolve together according to the rules of the basic coupling. Then*

$$\mathbf{P}^{\underline{\omega}}\{\underline{\omega}(t) = \underline{\zeta}(t)\} > 0 \qquad \text{for all } t > 0.$$

PROOF. The statement is that a neighboring second class antiparticle and second class particle will annihilate each other by time $t$ with positive probability. Intuitively this is obvious because one way for the annihilation to happen is that a brick is placed on column $i$ of $\underline{\omega}$ before any other jump in the vicinity of site $i$.

To make this rigorous, first note that the event $\{\underline{\omega}(t) = \underline{\zeta}(t)\}$ is increasing in time. Let $\underline{d}$ denote the height configuration of the initial $\underline{\zeta}$ configuration, normalized so that $d_0 = 0$. The height processes of $\underline{\omega}(t)$ and $\underline{\zeta}(t)$ are denoted by $\underline{h}(t)$ and $\underline{g}(t)$ as usual, with initial states $(\underline{h}, \underline{g}) = (\underline{h}(0), \underline{g}(0))$, but the distinct symbol $\underline{d}$ makes this next definition understandable. Define a bounded cylinder function $\varphi_{\underline{d}}$ on pairs of height configurations $(\underline{h}, \underline{g})$ by

$$\varphi_{\underline{d}}(\underline{h}, \underline{g}) := \mathbf{1}\{h_j = g_j = d_j \text{ for } i - 1 \leq j \leq i + 1\}.$$

Observe these two properties that follow from assumption (7.4):

- the function vanishes at time zero: $\varphi_{\underline{d}}(\underline{h}, \underline{g}) = 0$,
- at time zero the coupled generator applied to this function is positive:

$$(S(0)L\varphi_{\underline{d}})(\underline{h}, \underline{g}) = (L\varphi_{\underline{d}})(\underline{h}, \underline{g}) = r(\omega_i) - r(\zeta_i) + r(-\omega_{i+1}) - r(-\zeta_{i+1}) > 0$$

by (7.4) and *strict* monotonicity of the rates $r$.

By Lemma 6.2, $t \mapsto (S(t)L\varphi_{\underline{d}})(\underline{h}, \underline{g})$ is continuous up to some positive $T^{\underline{\omega}}$. (We are actually using the version of Lemma 6.2 that could be proved for a pair of coupled processes; see Remark 6.3.) Hence, by Theorem 2.6,

$$\mathbf{P}^{\underline{\omega}}\{h_j(t) = g_j(t) = g_j(0) \text{ for } i - 1 \leq j \leq i + 1\}$$
(7.5)
$$= \mathbf{E}^{\underline{\omega}}\varphi_{\underline{d}}(\underline{h}(t), \underline{g}(t))$$
$$= (S(t)\varphi_{\underline{d}})(\underline{h}, \underline{g}) = \int_0^t (S(s)L\varphi_{\underline{d}})(\underline{h}, \underline{g})\,ds > 0$$

for some time $t < T^{\underline{\omega}}$. The event in (7.5) implies that the second class particle and antiparticle have annihilated each other because among the columns $\{h_j, g_j : j = i - 1, i, i + 1\}$ it permits only one jump in column $h_i$ and no other move during time $[0, t]$. □

We are ready to prove ergodicity. According to Proposition 2.1 in [14], ergodicity is equivalent to the property that harmonic functions in $L^2(\underline{\mu}^{(\theta)})$ are a.e. constant. See also Corollaries 2 and 5 in IV. 2 of [12] for the same discussion in a discrete time setting. The next lemma completes the proof of Theorem 2.8.

LEMMA 7.3. *Let $\psi \in L^2(\underline{\mu}^{(\theta)})$ satisfy $S(t)\psi = \psi$ for all $t \geq 0$. Then $\psi$ is $\underline{\mu}^{(\theta)}$-a.s. constant.*

PROOF. $M(t) := \psi(\underline{\omega}(t))$ is an $L^2$-martingale. By stationarity,

$$0 = \mathbf{E}(M(t)^2) - \mathbf{E}(M(0)^2) = \mathbf{E}[(M(t) - M(0))^2] = \mathbf{E}[(\psi(\underline{\omega}(t)) - \psi(\underline{\omega}(0)))^2].$$

Consequently, $\psi(\underline{\omega}(t)) = \psi(\underline{\omega}(0))$ a.s. This we can express as

$$0 = \mathbf{E}[|\psi(\underline{\omega}(t)) - \psi(\underline{\omega}(0))|]$$
$$= \sum_{y,z \in \mathbb{Z}} \mathbf{E}[|\psi(\underline{\omega}(t)) - \psi(\underline{\omega}(0))||\omega_i(0) = y, \omega_{i+1}(0) = z]\mu^{(\theta)}(y)\mu^{(\theta)}(z).$$

Hence, for any $y, z \in \mathbb{Z}$,

(7.6) $\qquad \mathbf{E}[|\psi(\underline{\omega}(t)) - \psi(\underline{\omega}(0))||\omega_i(0) = y, \omega_{i+1}(0) = z] = 0.$

Let $\underline{\zeta}(t)$ be another process, initially coupled to $\underline{\omega}(0)$ according to (7.4). Let the pair $(\underline{\omega}(t), \underline{\zeta}(t))$ evolve according to the rules of the basic coupling. Then by the product structure of the initial measure and by (7.6),

$$\mathbf{E}[|\psi(\underline{\zeta}(t)) - \psi(\underline{\zeta}(0))||\zeta_i(0) = y - 1, \zeta_{i+1}(0) = z + 1]$$
$$= \mathbf{E}[|\psi(\underline{\omega}(t)) - \psi(\underline{\omega}(0))||\omega_i(0) = y - 1, \omega_{i+1}(0) = z + 1] = 0,$$

which implies

$$\mathbf{E}[|\psi(\underline{\zeta}(t)) - \psi(\underline{\zeta}(0))|]$$
$$= \sum_{y,z \in \mathbb{Z}} \mathbf{E}[|\psi(\underline{\zeta}(t)) - \psi(\underline{\zeta}(0))||\zeta_i(0) = y - 1,$$
$$\zeta_{i+1}(0) = z + 1]\mu^{(\theta)}(y)\mu^{(\theta)}(z) = 0.$$

Thus, remembering (7.4), we also have $\psi(\underline{\zeta}(t)) = \psi(\underline{\zeta}(0)) = \psi(\underline{\omega}(0)^{(i,i+1)})$ a.s. And now follows invariance of $\psi$ under deposition of a brick:

$$\int |\psi(\underline{\omega}^{(i,i+1)}) - \psi(\underline{\omega})| \cdot \mathbf{P}^{\underline{\omega}}\{\underline{\omega}(t) = \underline{\zeta}(t)\} \, d\underline{\mu}^{(\theta)}(\underline{\omega})$$
$$= \int \mathbf{E}^{\underline{\omega}}[|\psi(\underline{\zeta}(0)) - \psi(\underline{\omega}(0))| \cdot \mathbf{1}\{\underline{\omega}(t) = \underline{\zeta}(t)\}] \, d\underline{\mu}^{(\theta)}(\underline{\omega})$$
$$= \mathbf{E}[|\psi(\underline{\zeta}(t)) - \psi(\underline{\omega}(t))| \cdot \mathbf{1}\{\underline{\omega}(t) = \underline{\zeta}(t)\}] = 0.$$

By the previous lemma, $\mathbf{P}^{\underline{\omega}}\{\underline{\omega}(t) = \underline{\zeta}(t)\} > 0$ for all $\underline{\omega} \in \widetilde{\Omega}$, and so $\psi(\underline{\omega}^{(i,i+1)}) = \psi(\underline{\omega})$ $\underline{\mu}^{(\theta)}$-a.s.

We have shown that $\psi(\underline{\omega}^{(i,i+1)}) = \psi(\underline{\omega})$ $\underline{\mu}^{(\theta)}$-a.s. for all $i$. But invariance of $\psi$ under deposition of bricks implies invariance of $\psi$ under removal of bricks. By adding or removing bricks in the column $h_j$, we can interchange $\omega_j$ and $\omega_{j+1}$. Thus, $\psi$ is invariant under finite permutations, and then $\underline{\mu}^{(\theta)}$-a.s. constant by the Hewitt–Savage 0–1 law. □

## APPENDIX A: THE MEASURE $\mu^{(\theta)}$

We establish here some properties of $\mu^{(\theta)}$. Assumption (2.2) implies that

(A.1) $$\mathbf{E}^{(\theta)} e^{C|z|} < \infty \quad \text{for any constant } C.$$

In particular, $\mathbf{E}^{(\theta)} r(|z|)^p$ is finite for any $1 \leq p < \infty$.

LEMMA A.1. $\mathbf{E}^{(\theta)}(z)$ *is a strictly increasing function of $\theta$, and $\mathbf{E}^{(\theta)}(z) \to \pm\infty$ as $\theta \to \pm\infty$.*

REMARK A.2. $\mathbf{E}^{(\theta)}(z)$ represents particle density if particles and antiparticles are counted with opposite signs. By the lemma, the family $\{\underline{\mu}^{(\theta)}\}_{\theta \in \mathbb{R}}$ could also be indexed by density, rather than by the parameter $\theta$.

REMARK A.3. In the zero range process occupation numbers are nonnegative. Consequently, $\mathbf{E}^{(\theta)}(z) \geq 0$ and the lemma is valid only for $\theta \to \infty$. In the proof below consider only $\theta \geq 0$ and restrict the sums to $z \geq 0$.

PROOF OF LEMMA A.1. First, by (2.3),

$$Z(\theta) = \sum_{z=-\infty}^{\infty} \frac{e^{\theta z}}{r(|z|)!} = 1 + \sum_{z=1}^{\infty} \frac{e^{\theta z} + e^{-\theta z}}{r(z)!} \geq \sum_{z=1}^{\infty} \frac{e^{|\theta| z}}{r(z)!} \geq \sum_{z=1}^{\infty} \frac{e^{|\theta| z}}{\prod_{y=1}^{z} e^{\beta y}}$$

$$= \sum_{z=1}^{\infty} e^{(1/2)(2|\theta| z - \beta z^2 - \beta z)} = e^{(\beta/2)(|\theta|/\beta - 1/2)^2} \sum_{z=1}^{\infty} e^{-(\beta/2)(z - |\theta|/\beta + 1/2)^2}.$$

The sum is uniformly bounded from below by $e^{-9\beta/8}$ as $z$ will eventually get close to $|\theta|/\beta$. Now we proceed by

$$\frac{d}{d\theta} \ln(Z(\theta)) = \frac{\sum_{z=-\infty}^{\infty} z e^{\theta z}/(r(|z|)!)}{Z(\theta)} = \mathbf{E}^{(\theta)}(z),$$

$$\frac{d^2}{d\theta^2} \ln(Z(\theta)) = \frac{\sum_{z=-\infty}^{\infty} z^2 e^{\theta z}/(r(|z|)!) Z(\theta) - (\sum_{z=-\infty}^{\infty} z e^{\theta z}/(r(|z|)!))^2}{Z(\theta)^2}$$

$$= \mathbf{Var}^{(\theta)}(z) > 0.$$

The last line shows that the derivative $\mathbf{E}^{(\theta)}(z)$ of $\ln Z(\theta)$ is strictly increasing. Since $\ln Z(\theta)$ is bounded below by the parabola $\frac{1}{2}\beta(|\theta|/\beta - \frac{1}{2})^2 - 9\beta/8$, the derivative $\mathbf{E}^{(\theta)}(z)$ cannot be bounded either above or below. □

LEMMA A.2. *For any $\theta_1 \leq \theta_2$, there is a coupling measure $\underline{\mu}^{(\theta_1, \theta_2)}$ on $\Omega^2$ of which the first marginal is $\underline{\mu}^{(\theta_1)}$, the second marginal is $\underline{\mu}^{(\theta_2)}$, and*

$$\underline{\mu}^{(\theta_1, \theta_2)}\{(\underline{\omega}, \underline{\zeta}) : \omega_i \leq \zeta_i \text{ for all } i \in \mathbb{Z}\} = 1.$$

TABLE 5
*Rates for bricklayers at sites $\ell \leq i < \mathfrak{r}$ to lay brick to their right*

| With rate | $g_i$ | $f_i$ | $d_i$ | $d_{i+1}$ |
|---|---|---|---|---|
| $r(\xi_i) - r(\zeta_i)$ | | ↑ | ↓ | ↑ |
| $r(\zeta_i)$ | ↑ | ↑ | | |

PROOF. By the product structure, the lemma follows from the corresponding statement about the site-marginals, that is, about $\mu^{(\theta_1)}$ and $\mu^{(\theta_2)}$. It is enough to prove that for any nondecreasing, nonnegative function $f$ on $\mathbb{Z}$ (for which the expectations exist), we have

(A.2) $$\mathbf{E}^{(\theta_1)} f(z) \leq \mathbf{E}^{(\theta_2)} f(z),$$

see, for example, [8]. (A.2) is equivalent to showing that $\mathbf{E}^{(\theta)} f$ is nondecreasing in $\theta$. To this order, we consider

$$\frac{d}{d\theta} \mathbf{E}^{(\theta)} f(z) = \frac{d}{d\theta} \sum_{z=-\infty}^{\infty} f(z) \frac{1}{Z(\theta)} \cdot \frac{e^{\theta z}}{r(|z|)!}$$

$$= \sum_{z=-\infty}^{\infty} z \cdot f(z) \frac{1}{Z(\theta)} \cdot \frac{e^{\theta z}}{r(|z|)!}$$

$$- \sum_{z=-\infty}^{\infty} f(z) \frac{1}{Z(\theta)^2} \cdot \frac{e^{\theta z}}{r(|z|)!} \sum_{y=-\infty}^{\infty} y \cdot \frac{e^{\theta y}}{r(|y|)!}$$

$$= \mathbf{E}^{(\theta)}(z \cdot f(z)) - \mathbf{E}^{(\theta)} f(z) \cdot \mathbf{E}^{(\theta)} z.$$

The previous display is nonnegative because $f$ is nondecreasing, hence, $z$ and $f(z)$ are positively correlated. □

## APPENDIX B: THE BASIC COUPLING

We describe the basic coupling of the finite-volume $(\ell, \mathfrak{r}, \theta)$-processes in this section.

TABLE 6
*Rates for bricklayers at sites $\ell < i \leq \mathfrak{r}$ to lay brick to their left*

| With rate | $g_{i-1}$ | $f_{i-1}$ | $d_{i-1}$ | $d_i$ |
|---|---|---|---|---|
| $r(-\zeta_i) - r(-\xi_i)$ | ↑ | | ↑ | ↓ |
| $r(-\xi_i)$ | ↑ | ↑ | | |

TABLE 7
*Rates for the bricklayer at site $\mathfrak{r}$ to lay brick to his right*

| With rate | $g_{\mathfrak{r}}$ | $f_{\mathfrak{r}}$ | $d_{\mathfrak{r}}$ | $d_{\mathfrak{r}+1}$ |
|---|---|---|---|---|
| $r(\xi_{\mathfrak{r}}) - r(\zeta_{\mathfrak{r}})$ | | ↑ | ↓ | ↑ |
| $e^{-\theta_1} - e^{-\theta_2}$ | ↑ | | ↑ | ↓ |
| $r(\zeta_{\mathfrak{r}}) + e^{-\theta_2}$ | ↑ | ↑ | | |

THE COUPLINGS FOR LEMMA 4.1. We show the coupling which preserves $\xi_i(t) \geq \zeta_i(t)$ on sites $\ell \leq i \leq \mathfrak{r}$ for all times $t > 0$. We denote the height of $\underline{\zeta}$ and $\underline{\xi}$ by $g$ and $f$, respectively; the number of second class particles is $d_i := \xi_i - \zeta_i$. We do not have antiparticles when starting the processes. We rewrite Tables 3 and 4 for inner sites to Tables 5 and 6 with the present notation. For sites $\ell$ and $\mathfrak{r}$, where we modified the rates, the processes are coupled according to Tables 7 and 8.

These tables are valid while $\xi_i \geq \zeta_i$ holds for all $\ell \leq i \leq \mathfrak{r}$. However, they preserve this condition as no antiparticles are created according to these tables. Decrease of $d_i$ can only happen where $\xi_i > \zeta_i$, that is, for sites where there is particle to jump from.

Now, by the same method, we couple $\underline{\zeta}$ and $\underline{\eta}$, where the initial distribution of $\underline{\eta}$ is $\underline{\mu}^{(\theta_1)}$, and both processes evolve according to the $(\ell, \mathfrak{r}, \theta_1)$-generator. Writing $\eta$ instead of $\zeta$, $\zeta$ instead of $\xi$ and $\theta_1$ in place of $\theta_2$ as well makes us possible to repeat our arguments and to conclude $\zeta_i(t) \geq \eta_i(t)$ for all sites $\ell \leq i \leq \mathfrak{r}$ and $t > 0$. Hence, we see that $\eta_i(t) \leq \zeta_i(t) \leq \xi_i(t)$, where $\eta_i(t)$ and $\xi_i(t)$ have distributions $\mu^{(\theta_1)}$ and $\mu^{(\theta_2)}$, respectively, as these are processes started and evolving in their stationary distributions. □

**Acknowledgments.** The authors thank Maury Bramson, József Fritz, Thomas G. Kurtz, Ellen Saada, Bálint Tóth and S. R. Srinivasa Varadhan for fruitful conversations and advice. We also wish to thank the referee for careful readings and valuable suggestions about the paper. A large part of this work was done while Firas Rassoul–Agha was at the Mathematical

TABLE 8
*Rates for the bricklayer at site $\ell$ to lay brick to his left*

| With rate | $g_{\ell-1}$ | $f_{\ell-1}$ | $d_{\ell-1}$ | $d_\ell$ |
|---|---|---|---|---|
| $r(-\zeta_\ell) - r(-\xi_\ell)$ | ↑ | | ↑ | ↓ |
| $e^{\theta_2} - e^{\theta_1}$ | | ↑ | ↓ | ↑ |
| $r(-\xi_\ell) + e^{\theta_1}$ | ↑ | ↑ | | |

Biosciences Institute, Ohio State University. A preliminary version of Sections 1–3, only for the zero range process and with partial proofs, appeared in [4]. The present paper delivers the complete proofs and further results anticipated in [4].

M. BALÁZS
DEPARTMENT OF STOCHASTICS
INSTITUTE OF MATHEMATICS
BUDAPEST UNIVERSITY
OF TECHNOLOGY AND ECONOMICS
1 EGRY JÓZSEF U.
1111 BUDAPEST
HUNGARY
E-MAIL: balazs@math.bme.hu

F. RASSOUL-AGHA
DEPARTMENT OF MATHEMATICS
UNIVERSITY OF UTAH
155 SOUTH 1400 EAST
SALT LAKE CITY, UTAH 84112-0090
USA
E-MAIL: firas@math.utah.edu

T. SEPPÄLÄINEN
MATHEMATICS DEPARTMENT
UNIVERSITY OF WISCONSIN–MADISON
VAN VLECK HALL 480 LINCOLN DR
MADISON, WISCONSIN 53706-1388
USA
E-MAIL: seppalai@math.wisc.edu

S. SETHURAMAN
DEPARTMENT OF MATHEMATICS
IOWA STATE UNIVERSITY
430 CARVER
AMES, IOWA 50011
USA
E-MAIL: sethuram@iastate.edu